\documentclass{article}

\usepackage{color}
\usepackage{graphicx, graphics}
\usepackage{subfigure}
\usepackage{amsmath,amstext}  
\usepackage{amsfonts,amssymb, bbm} 
\usepackage[LGR, T1]{fontenc}
\usepackage[utf8]{inputenc}
\usepackage[english]{babel}
\usepackage{epsfig}
\usepackage{epsf}
\usepackage{color}
\usepackage{vmargin}

\definecolor{dkg}{rgb}{0,0.5,0}

\numberwithin{equation}{section}

\newcommand{\dN}{\mathbb{N}}
\newcommand{\dB}{\mathbb{B}}

\newcommand{\dT}{\mathbb{T}}
\newcommand{\dG}{\mathbb{G}}
\newcommand{\dF}{\mathbb{F}}
\newcommand{\dR}{\mathbb{R}}
\newcommand{\dE}{\mathbb{E}}
\newcommand{\dP}{\mathbb{P}}
\newcommand{\cA}{\mathcal{A}}
\newcommand{\cN}{\mathcal{N}}

\newcommand{\cB}{\mathcal{B}}

\newcommand{\cO}{\mathcal{O}}
\newcommand{\cR}{\mathcal{R}}
\newcommand{\cT}{\mathcal{T}}

\newcommand{\cF}{\mathcal{F}}
\newcommand{\cE}{\mathcal{E}}
\newcommand{\cV}{\mathcal{V}}
\newcommand{\cW}{\mathcal{W}}

\newcommand{\rI}{\boldsymbol{\mathrm{I}}}

\newcommand{\veps}{\varepsilon}
\newcommand{\eps}{\epsilon}
\newcommand{\wh}{\widehat}

\newcommand{\inds}[1]{\mathbbm{1}_{#1}}
\newcommand{\ind}[1]{\mathbbm{1}_{\{#1\}}}
\newcommand{\indnex}{\mathbbm{1}_{\overline{\cE}}}
\newcommand{\bs}[1]{\boldsymbol{#1}}
\newcommand{\cM}{<\!\!\bs{M}\!\!>}
\newcommand{\bX}{X^*}

\font\calcal=cmsy10 scaled\magstep1
\def\build#1_#2^#3{\mathrel{\mathop{\kern 0pt#1}\limits_{#2}^{#3}}}
\def\liml{\build{\longrightarrow}_{}^{{\mbox{\calcal L}}}}

\def\videbox{\mathbin{\vbox{\hrule\hbox{\vrule height1ex \kern.5em
\vrule height1ex}\hrule}}}

\newtheorem{defi}{Definition}[section]
\newtheorem{Lemma}[defi]{Lemma}
\newtheorem{Corollary}[defi]{Corollary}
\newtheorem{Proposition}[defi]{Proposition}

\newtheorem{Theorem}[defi]{Theorem}
\title{Random coefficients bifurcating autoregressive processes}
\author{ Beno\^{\i}te de Saporta\\
Univ. Bordeaux, Gretha, UMR 5113, IMB, UMR 5251, F-33400 Talence, France\\
CNRS, Gretha, UMR 5113, IMB, UMR 5251, F-33400 Talence, France\\
INRIA Bordeaux Sud Ouest, team CQFD, F-33400 Talence, France\\
\and Anne G\'egout-Petit\\
Univ. Bordeaux, IMB, UMR 5251, F-33400 Talence, France\\
CNRS, IMB, UMR 5251, F-33400 Talence, France\\
INRIA Bordeaux Sud Ouest, team CQFD, F-33400 Talence, France\\
\and Laurence Marsalle \\
Univ. Lille 1, Laboratoire Paul Painlev\'e, UMR~8524, F-59 655 Villeneuve d'Ascq, France\\
CNRS, Laboratoire Paul Painlev\'e, UMR~8524, F-59 655 Villeneuve d'Ascq, France}

\begin{document}

\maketitle

\begin{abstract}
This paper presents a new model of asymmetric bifurcating autoregressive process with random coefficients. We couple this model with a Galton-Watson tree to take into account possibly missing observations. We propose least-squares estimators for the various parameters of the model and prove their consistency, with a convergence rate, and asymptotic normality. We use both the bifurcating Markov chain and martingale approaches and derive new results in both these frameworks.
\end{abstract}

\section{Introduction}
\label{section intro}
In the 80's, Cowan and Staudte \cite{CoSt86} introduced Bifurcating Autoregressive processes (BAR) as a parametric model to study cell lineage data. A quantitative characteristic of the cells (e.g. growth rate, age at division) is recorded over several generations descended from an initial cell, keeping track of the genealogy to study inherited effects. As a cell usually gives birth to two offspring by division, such genealogies are naturally structured as binary trees. BAR processes are thus a generalization of autoregressive processes (AR) to this binary tree structure, by modeling each line of descent as a first order AR process, allowing the environmental effects on sister cells to be correlated. Statistical inference for the parameters of BAR processes has been widely studied, either based on the observation of a single tree growing to infinity \cite{CoSt86, Hug96, HuBa00, ZhBa05b} or on a large number of small independent trees \cite{HuSt94, HuBa99}. See also \cite{HwBa11, HwBa09} for processes indexed by general trees.

Various extensions of the original model have been proposed, e.g. non gaussian noise sequence \cite{BaZh04, ZhBa05a}, higher order AR \cite{HuBa00, ZhBa05a} or moving average AR \cite{HuBa99}. Since 2005, evidence of asymmetry in cell division has been established by biologists \cite{SMPT05} and an asymmetric BAR model has been introduced by Guyon \cite{Guy07} where the coefficients of the AR processes of sister cells are allowed to be different. This model was further extended to higher order AR \cite{BSG09}, to take missing data into account \cite{DM08, SGM11, SPL12} and with parasite infection \cite{Ban08}. 

To the best of our knowledge, only two papers \cite{BH99} and \cite{Vassili12} deal with random coefficient BAR processes. In the former by Bui and Huggins it is explained that random coefficients BAR processes can account for observations that do not fit the usual BAR model. For instance, the extra randomness can model irregularities in nutrient concentrations in the media in which the cells are grown. Other evidence for the need of richer models can be found e.g. in \cite{GBPSDT05}. In this paper, we propose a new model for random coefficient BAR processes (R-BAR). It is more general than that of Bui and Huggins, as the random variables are not supposed to be Gaussian, they may not have moments of all order and correlation between all the sources of randomness are allowed. Moreover, we propose an asymmetric model in the continuation of \cite{Guy07, BSG09, DM08, SGM11, SPL12, Vassili12} in the context of missing data. Indeed, experimental data are often incomplete and it is important to take this phenomenon into account for the inference. As in \cite{DM08, SGM11} we model the structure of available data by a Galton-Watson tree, instead of a complete binary tree. Our model is close to that developed in \cite{Vassili12}, but the assumptions on the noise process are different as we allow correlation between the two sources of randomness but require higher moments because of the missing data and because we do not use a weighted estimator. The main difference is that the model in \cite{Vassili12} is fully observed, whereas ours allows for missing observations.

Our approach for the inference of our model is also different from \cite{BH99, Vassili12}. As we cannot use maximum likelihood estimation, we propose modified least squares estimators as in \cite{NiQu82}. In \cite{BH99}, inference is based on an asymptotically infinite number of small replicated trees. Here, as in \cite{Vassili12}, we consider one single tree growing to infinity but our least squares estimator is not weighted. The originality of our approach is that it combines the bifurcating Markov chain and martingale approaches. Bifurcating Markov chains (BMC) were introduced in \cite{Guy07} on complete binary trees and further developed in \cite{DM08} in the context of missing data on Galton-Watson trees. BAR models can be seen as a special case of BMC. This interpretation allows us to establish the convergence of our estimators. A by-product of our procedure is a new general result for BMC on Galton-Watson trees. Indeed, in \cite{Guy07, DM08} the driven noise sequence is assumed to have moments of all order. Here, we establish new laws of large numbers for polynomial functions of the BMC where the noise sequence only has moments up to a given order. The strong law of large numbers \cite{Duflo97} and the central limit theorem \cite{HaHe80, Ham94, Duflo97} for martingales have been previously used in the context of BAR processes \cite{BaZh04, ZhBa05a, ZhBa05b} and adapted to special cases of martingales on binary trees \cite{BSG09, SGM11, SPL12, Vassili12}. In this paper, we establish a general law of large numbers for square-integrable martingales on Galton-Watson binary trees. This result is applied to our R-BAR model to obtain sharp convergence rates and a quadratic strong law for our estimators.

The paper is organized as follows. In Section~\ref{section model}, we give the precise definition of our R-BAR model on a Galton-Watson tree and state our main assumptions. In Section~\ref{section estimation}, we give modified least squares estimators and state the convergence results we obtained: consistency with convergence rate and asymptotic normality. In Section~\ref{section BMC}, we recall the BMC framework, prove a new law of large numbers under limited moment conditions and apply it to our R-BAR model to derive the consistency of our estimators. In Section~\ref{section martingales} we establish a new general law of large numbers for square-integrable martingales on Galton-Watson trees and use it to derive convergence rates and quadratic strong laws for our estimators. In Section~\ref {section TCL} we establish the asymptotic normality by using central limit theorems for martingales. Finally in Section~\ref{section data} we apply our estimation procedure to the E. coli data of \cite{SMPT05}.
\section{Model}
\label{section model}
In the sequel, all random variables are defined on the probability state space $(\Omega, \cA, \dP)$. As in the previous literature, we use the index 1 for the original cell, and the two offspring of cell $k$ are labelled $2k$ and $2k+1$. Consider the first-order asymmetric random coefficients bifurcating autoregressive process (R-BAR) given, for all $k\geq 1$, by
\begin{equation}\label{defbar}
\left\{
    \begin{array}{lcrcrcrcr}
     X_{2k} & = & (b_{2k}&+& \eta_{2n})X_k &+ &(a_{2k}&+ &\varepsilon_{2k}), \\
     X_{2k+1} & = &  (b_{2k+1}&+& \eta_{2k+1})X_k &+ & (a_{2k+1}&+&\varepsilon_{2k+1}),
    \end{array}\right.
\end{equation}
with the following notations: for all $k\geq 1$, 
\begin{equation*}
\left\{
    \begin{array}{ll}
     a_{2k} = a, \vspace{1ex}\\
     b_{2k} = b,
    \end{array}\right.
\hspace{1cm}
\text{and}
\hspace{1cm}
\left\{
    \begin{array}{ll}
     a_{2k+1} = c, \vspace{1ex}\\
     b_{2k+1} = d.
    \end{array}\right.
\end{equation*}
The initial state $X_1$ is the characteristic of the original ancestor while the sequence $(\veps_{2k},\eta_{2k},\veps_{2k+1},\eta_{2k+1})_{n\geq 1}$
is the driving noise of the process, and the parameter $(a,b,c,d)$ belongs to $\dR^4$.
One can see this
R-BAR process as a random-coefficient first-order autoregressive process on a binary tree,
where each vertex represents an individual or cell, vertex $1$ being the original ancestor.
For all $n\geq 1$, denote the n-th generation by
$$\dG_n = \{2^n, 2^n+1,\ldots,2^{n+1}-1\}.$$
In particular, $\dG_0 = \{1\}$ is the initial generation and $\dG_1 = \{2,3\}$ is
the first generation of offspring from the original ancestor. 
Finally, denote by
$$\dT_n = \bigcup_{\ell=0}^n\dG_{\ell},$$
the sub-tree of all individuals from the original individual up to the $n$-th generation and $\dT$ the complete tree. Note that the cardinality $|\dG_n|$ of $\dG_n$  is $2^n$ while that
of $\dT_n$ is $|\dT_n|=2^{n+1}-1$.
In all the sequel, we shall use the  following  hypotheses.
\begin{description}
\item[(H.1)] The sequence $(\veps_{2k},\eta_{2k},\veps_{2k+1},\eta_{2k+1})_{k\geq 1}$ is independent and identically distributed. It is also independent from  $X_1$.
\item[(H.2)] The random variables $\veps_2$, $\eta_2$, $\veps_3$, $\eta_3$ and $X_1$ have moments of all order up to $4\gamma$, for some $\gamma\geq1$. The following hypotheses will be used
\begin{eqnarray*}
&\dE[\veps_{2}]  = \dE[\veps_{3}]  = 0,
\dE[\veps_{2}^2]= \dE[\veps_{3}^2]=\sigma_{\veps}^2>0
\hspace{0.25cm}\text{and}\hspace{0.25cm}
\dE[\veps_{2}\veps_{3}]= \rho_{\veps},&\\
&\dE[\eta_{2}] = \dE[\eta_{3}] =0,
\dE[\eta_{2}^2]=\dE[\eta_{3}^2]= \sigma_{\eta}^2>0
\hspace{0.25cm}\text{and}\hspace{0.25cm}
\dE[\eta_{2}\eta_{3}]= \rho_{\eta},&\\
&\dE[\veps_{2+i}\eta_{2+j}]= \rho_{ij}, \textrm{\ for\ } (i,j)\in\{0,1\},
\hspace{0.15cm}\text{and}\hspace{0.15cm}
\rho = \frac{1}{2}(\rho_{01}+\rho_{10}).&
\end{eqnarray*}
\end{description}
When dealing with the biological issue of cell lineages, it may happen that a lineage is incomplete. Indeed, cells may die or measurements may be impossible or faulty on some cells. Taking into account such a phenomenon, we introduce the observation process, $(\delta_k)_{k \in \dT}$. 
We use the same framework as in \cite{DM08}, and not the more general one introduced in \cite{SGM11}.
Basically, $\delta_k=1$ if cell $k$ is observed, $\delta_k=0$ otherwise.  We set $\delta_1=1$ and define the whole sequence through the following equalities:
\begin{equation}
 \delta_{2k}=\delta_k \xi_k^0 \quad \text{and} \quad \delta_{2k+1}=\delta_k \xi_k^1, 
\end{equation}
where the sequence $\big(\bs{\xi}_k=(\xi_k^0, \xi_k^1)\big)_{k \in \dT}$ is a sequence of independent identically distributed random vectors of $\{0,1\}^2$ whose common distribution is specified by the following generating function
$$\dE\big[s_0^{\xi_1^0}s_1^{\xi_1^1}\big]= (1-p_0-p_1-p_{01}) + p_0 s_0 + p_1 s_1 + p_{01} s_0s_1.$$
We also suppose that the observation process is independent from the state process $(X_n)$.
\begin{description}
 \item [(H.3)] The sequence $(\bs{\xi}_k)_{k\in\dT}$ is independent from $(\veps_{2k},\eta_{2k},\veps_{2k+1},\eta_{2k+1})_{k\in\dT}$ and from $X_1$.
\end{description}
Notice that the process $(\delta_k)_{k \in \dT}$ takes its values in $\{0,1\}$, and that if $k \in \dT$ is such that $\delta_k=0$, then $\delta_{2^nk+i} = 0$, for all $i \in \{0, \ldots, 2^n -1\}$ and all $n \ge 1$. So to speak, if individual $k$ is not observed, all its descendants are also missing. We now define the sets of observed data
$$\dG_n^*  = \{k \in \dG_n : \delta_k=1\} \quad \text{and} \quad \dT_n^* = \{k \in \dT_n : \delta_k=1\} = \cup_{\ell=0}^n \dG_{\ell}^*.$$
Thanks to the i.i.d. property of the sequence $(\bs{\xi}_k)$, the sequence of cardinalities $(|\dG_n^*|)_{n\ge 0}$ is a Galton-Watson (GW) process with reproduction generating function
$$z \mapsto (1-p_0-p_1-p_{01})  + (p_0+p_1)z + p_{01}z^2,$$
and mean $m = 2p_{01}+ p_0+p_1$. The following equalities thus hold (see e.g. \cite{Har63})
$$\dE[|\dG_n^*|]= m ^n \quad \text{and} \quad \dE[|\dT_n^*|] = \sum_{\ell=0}^n \dE[|\dG_{\ell}^*|] = \frac{m^{n+1}-1}{m-1}\cdot$$
According to the position of the mean $m$ of the reproduction law with respect to $1$, it is well known that the population becomes extinct or not. More precisely, if $m \le 1$ then we have extinction almost surely, in the sense that $\dP(\cup_{n \ge 0} \{ |\dG_n^*|=0\}) = 1$. But if $m > 1$, there is a positive probability of survival of the population: $\dP(\cap_{n \ge 0} \{ |\dG_n^*|>0\}) >0$. This latter case is called the super-critical case, and we assume that we are in that case.
\begin{description}
 \item[(H.4)] The mean of the reproduction law is greater than 1: $m >1$. 
\end{description}
On the set of non-extinction, the growth of the population is exponential, more precisely there exists some non-negative square-integrable random variable $W$ such that
\begin{equation}
\label{eq cvGW}
 \lim_{n \rightarrow \infty} \frac{|\dG_n^*|}{m^n} = W \ \text{a.s.,}\qquad \text{and} \qquad \{W>0\} = \cap_{n \ge 0} \{ |\dG_n^*|>0\} \ \text{a.s.}
\end{equation}
This immediately entails that
\begin{equation} \label{eq cvGW2}
 \lim_{n \rightarrow \infty} \frac{|\dT_n^*|}{m^n} = W \times \frac{m}{m-1}\quad \text{a.s.}
\end{equation}
We will denote by $\cE$ the extinction set $\cE=\cup_{n \ge 0} \{ |\dG_n^*|=0\}$ and $\overline{\cE}$ its complementary set. Note that under assumption \textbf{(H.4)}, $\overline{\cE}$ has a positive probability: $\dP({\overline{\cE}})>0$.
We need one more assumption combining the R-BAR and GW processes. 
\begin{description}
 \item [(H.5)] There exist $1\leq\kappa\leq\gamma$ such that 
 \begin{equation*}
 \frac{p_0+p_{01}}{m}\mathbb{E}[(b+\eta_2)^{4\kappa}]+ \frac{p_1+p_{01}}{m}\mathbb{E}[(d+\eta_3)^{4\kappa}]<1.
 \end{equation*}
\end{description}
This is the analogous of the usual stability assumption for the autoregression expressed by $\max\{|b|,|d|\}<1$ in the case of deterministic coefficients. It ensures that the values of $|X_k|$ do not tend to infinity. Note that the assumption above is slightly weaker than the one for deterministic coefficients. Indeed, in the fully observed case and when $\eta_2=\eta_3=0$, this equation reduces to $(b^{4\kappa}+d^{4\kappa})/2<1$. The special form of this assumption is explained in Section~\ref{subsection ergo} and is closely linked to the properties of the R-BAR process as a bifurcating Markov chain.

Finally, denote by $\dF=(\cF_n)$ the natural filtration of the R-BAR process $(X_k)_{k\in\dT}$, which means that $\cF_n$ is the $\sigma$-algebra generated by all individuals up to the
$n$-th generation, $\cF_n = \sigma\{X_k, k\in \dT_n\}$. We also introduce the sigma field $\cO = \sigma \{\delta_k , k\in \dT \}$ generated by the observation process. We shall assume that all the history of the observation process $(\delta_k)$ is known at time $0$ and use the filtration $\dF^{\cO}=(\cF^{\cO}_n)$ defined for all $n$ by
\begin{equation*}
\cF^{\cO}_n=\cO\vee\sigma\{\delta_kX_k, k\in\dT_n\}=\cO\vee\sigma\{X_k, k\in\dT_n^*\}.
\end{equation*}
Note that $\cF^{\cO}_n$ is a sub-$\sigma$-field of $\cO\vee\cF_n$.
\section{Estimation}
\label{section estimation}
We now give some least-squares estimators of our  parameters and state our main results on their asymptotic behavior.
\subsection{Estimators}
\label{section estimators}
We propose to use the standard least-squares (LS) estimator $\wh{{\bs{\theta}}}_n=(\wh{a}_n,\wh{b}_n,\wh{c}_n,\wh{d}_n)^t$ of $\bs{\theta}= (a,b,c,d)^t$ which minimizes the following expression
\begin{equation*}
\Delta_n(\bs{\theta})=\frac{1}{2}\sum_{k\in \dT_{n-1}}\delta_{2k} (X_{2k}-a-bX_k)^2 + \delta_{2k+1}(X_{2k+1}-c-dX_k)^2.
\end{equation*}
Consequently, for all $n\geq 1$ the following equality holds
\begin{equation*}
\wh{{\bs{\theta}}}_n\ =\
\bs{S}_{n-1}^{-1}
\sum_{k \in \mathbb{T}_{n-1}}
\left( \begin{array}{c}
\delta_{2k} X_{2k}  \\
\delta_{2k} X_kX_{2k} \\
\delta_{2k+1} X_{2k+1} \\
\delta_{2k+1} X_kX_{2k+1}
\end{array}\right),\ \textrm{with}\quad
\bs{S}_{n-1} =  \left( \begin{array}{cc}
\bs{S}^0_{n-1} & 0 \\
0 & \bs{S}^1_{n-1}
\end{array} \right),
\end{equation*}
and
$
\bs{S}^0_{n-1} =\sum_{k \in \mathbb{T}_{n-1}}\delta_{2k}\left(
\begin{array}{cc}
1 & X_k \\
X_k &X^2_k
\end{array}\right)$, 
$\bs{S}^1_{n-1} =\sum_{k \in \mathbb{T}_{n-1}}\delta_{2k+1}\left( \begin{array}{cc}
1 & X_k \\
X_k & X^2_k
\end{array}\right)$.

We now turn to the estimation of the parameters of the conditional covariance of $(\varepsilon_{2},\eta_{2},\varepsilon_{3},\eta_{3})$. Following \cite{NiQu82}, we obtain a modified least squares estimator of 
$\bs{\sigma}=(\sigma^2_{\veps}, \rho_{00}, \rho_{11}, \sigma^2_{\eta})^t$ by minimizing
\begin{equation*}
\Delta'_n(\bs{\sigma})=\frac{1}{2}\sum_{\ell=1}^{n-1}\sum_{k\in \dG_{\ell}}( \wh{\eps}^2_{2k}-\dE[ {\eps}^{2}_{2k}|\cF_{\ell}^{\cO}])^2 + ( \wh{\eps}^2_{2k+1}-\dE[ {\eps}^{2}_{2k+1}|\cF_{\ell}^{\cO}])^2,
\end{equation*}
where for all $k\in\mathbb{G}_n$, 
\begin{equation*}
\left\{
    \begin{array}{lcl}
     {\eps}_{2k} &=&\delta_{2k} (\veps_{2k} + \eta_{2k}X_k), \vspace{1ex}\\
     {\eps}_{2k+1} &=&\delta_{2k+1} ( \veps_{2k+1} + \eta_{2k+1}X_k),
    \end{array}\right.
    \  \left\{
    \begin{array}{lcl}
     \wh{\eps}_{2k} &=&\delta_{2k} ( X_{2k} - \wh{a}_n - \wh{b}_{n}X_k), \vspace{1ex}\\
     \wh{\eps}_{2k+1} &=& \delta_{2k} (X_{2k+1} - \wh{c}_{n} - \wh{d}_{n}X_k).
    \end{array}\right.
\end{equation*}
Under assumptions \textbf{(H.2)} and \textbf{(H.3)}, one obtains the following estimator
\begin{equation}
\label{sigma}
\wh{{\bs{\sigma}}}_n\ =\ 
\bs{U}^{-1}_{n-1}
\sum_{k \in \mathbb{T}_{n-1}}
\left( 
\wh{\eps}^2_{2k} +\wh{\eps}^2_{2k+1},
2X_k\wh{\eps}^2_{2k} ,
2X_k \wh{\eps}^2_{2k+1},
X^2_k( \wh{\eps}^2_{2k} + \wh{\eps}^2_{2k+1} )
\right)^t,
\end{equation}
where
\begin{equation*}
\bs{U}_n = \sum_{k \in \mathbb{T}_{n}}\left( \begin{array}{cccc}
\delta_{2k}+\delta_{2k+1}  & 2\delta_{2k}X_k & 2\delta_{2k+1}X_k & (\delta_{2k}+\delta_{2k+1})X^2_k\\
2\delta_{2k}X_k  & 4\delta_{2k}X^2_k & 0 &2\delta_{2k} X^3_k \\
2\delta_{2k+1}X_k  & 0 & 4\delta_{2k+1}X^2_k & 2\delta_{2k+1}X^3_k\\
(\delta_{2k}+\delta_{2k+1)}X_k^2  & 2\delta_{2k}X^3_k & 2\delta_{2k+1}X^3_k  & (\delta_{2k}+\delta_{2k+1)}X^4_k
\end{array}\right).
\end{equation*}
Note that if $\sigma^2_{\eta}=0$ the estimator of $\sigma^2_{\veps}$ above corresponds to the empirical estimator already used in \cite{SGM11}.
Similarly, the least-squares estimator of $\bs{{\rho}}=(\rho_{\veps}, \rho, \rho_{\eta})^t$ minimizes
\begin{equation*}
\Delta''_n(\bs{{\rho}})=\frac{1}{2}\sum_{\ell=1}^{n-1}\sum_{k\in \dG_{\ell}}( \wh{\eps}_{2k} \wh{\eps}_{2k+1}-\dE[ {\eps}_{2k} {\eps}_{2k+1}|\cF_{\ell}^{\cO}])^2,
\end{equation*}
and one obtains
\begin{equation}
\label{rho}
\wh{\bs{\rho}}_n= 
\bs{V}^{-1}_{n-1}
\sum_{k \in \mathbb{T}_{n-1}}
\left(
\wh{\eps}_{2k} \wh{\eps}_{2k+1},
2X_k \wh{\eps}_{2k}  \wh{\eps}_{2k+1},
X^2_k\wh{\eps}_{2k}  \wh{\eps}_{2k+1} 
\right)^t,
\end{equation}
where
\begin{equation*}
\bs{V}_n = \sum_{k \in \mathbb{T}_{n}}\delta_{2k}\delta_{2k+1}\left( \begin{array}{cccc}
1 &2 X_k & X^2_k  \\
2X_k & 4X^2_k & 2X^3_k  \\
X_k^2 & 2X^3_k & X^4_k \\
\end{array}\right).
\end{equation*}
Note that one cannot identify $\rho_{01}$ from $\rho_{10}$, hence the use of $\rho=(\rho_{01}+\rho_{10})/2$.
Again if $\sigma^2_{\eta}=0$, we retrieve the empirical estimator of  $\rho_{\veps}$ used in  \cite{SGM11}.

\subsection{Main results}
\label{section main res}
We now state our main results. The first one establishes the consistency of our estimators on the non-extinction set.
\begin{Theorem}\label{th consist}
Under assumptions \emph{\textbf{(H.1-5)}}, and if $\kappa\geq 2$, the following convergence holds
\begin{equation*}
\lim_{n\rightarrow\infty}\ind{|\dG_n^*|>0}\bs{\wh{\theta}}_n=\bs{\theta}\indnex\quad \text{a.s.}
\end{equation*}
and if in addition $\kappa\geq 4$ then the following convergences also hold
\begin{equation*}
\lim_{n\rightarrow\infty}\ind{|\dG_n^*|>0}\bs{\wh{\sigma}}_n=\bs{\sigma}\indnex\ \text{a.s.},\qquad \lim_{n\rightarrow\infty}\ind{|\dG_n^*|>0}\bs{\wh{\rho}}_n=\bs{\rho}\indnex\ \text{a.s.}
\end{equation*}
\end{Theorem}
The next results give convergence rates for the estimators.
\begin{Theorem}\label{th cv theta chap}
Under assumptions \emph{\textbf{(H.1-5)}} and if $\kappa\geq 4$, for all $\delta>1/2$,  the following convergence rate holds
\begin{equation*}
\|\wh{\bs{\theta}}_n-\bs{\theta} \|^2=o(n^{\delta}m^{-n})\quad\textrm{a.s.}\
\end{equation*}
with the quadratic strong law
\begin{equation*}
\lim_{n\rightarrow\infty}\ind{|\dG_n^*|>0}\frac{1}{n}\sum_{\ell=1}^n|\dT_{\ell-1}^*|^{-1}(\wh{\bs{\theta}}_{\ell}-\bs{\theta})^t\bs{S}\bs{\Sigma}^{-1}\bs{S}(\wh{\bs{\theta}}_{\ell}-\bs{\theta})=tr(\bs{\Gamma}\bs{\Sigma}^{-1})\indnex\qquad\textrm{a.s.}
\end{equation*}
where $\bs{S}$, $\bs{\Gamma}$ and  $\bs{\Sigma}$ are $4\times4$ matrices defined respectively in Proposition~\ref{prop lim SUV}, Lemma~\ref{lem lim Gamma} and Lemma~\ref{lem lim Sigma}.
\end{Theorem} 
For all $n$, set
\begin{equation*}
{{\bs{\sigma}}}_n\ = 
\bs{U}^{-1}_{n-1}
\sum_{k \in \mathbb{T}_{n-1}}
\left( \begin{array}{c}
{\eps}^2_{2k} +{\eps}^2_{2k+1} \\
2X_k{\eps}^2_{2k} \\
2X_k{\eps}^2_{2k+1}\\
X^2_k( {\eps}^2_{2k} + {\eps}^2_{2k+1} )
\end{array}\right),\ 
{{\bs{\rho}}}_n\ = 
\bs{V}^{-1}_{n-1}
\sum_{k \in \mathbb{T}_{n-1}}
\left( \begin{array}{c}
{\eps}_{2k} {\eps}_{2k+1} \\
2X_k {\eps}_{2k}  {\eps}_{2k+1}\\
X^2_k{\eps}_{2k}  {\eps}_{2k+1} 
\end{array}\right).
\end{equation*}
\begin{Theorem}\label{th sigma}
Under assumptions \emph{\textbf{(H.1-5)}} and if $\kappa\geq 8$, the following convergences hold
\begin{equation*}
\lim_{n\rightarrow\infty} \ind{|\dG_n^*|>0}{{\bs{\sigma}}}_n={{\bs{\sigma}}}\indnex\qquad a.s.
\end{equation*}
and
\begin{eqnarray*}
{\lim_{n\rightarrow\infty}\ind{|\dG_n^*|>0}\frac{|\dT_{n-1}^*|}{n}(\wh{{\bs{\sigma}}}_n-{{\bs{\sigma}}}_n)}
&=&\bs{U}^{-1}\left(q_0(0)+q_1(0),\ 2q_0(1),\ 2q_1(1),\ q_0(2)+q_1(2)\right)^t\indnex\qquad a.s.
\end{eqnarray*}
where $\bs{U}$ is a $4\times4$ matrix defined in Proposition~\ref{prop lim SUV} and the $q_i(r)$ are scalars defined in Lemmas~\ref{lem eps2-0}, \ref{lem eps2-1}, \ref{lem esp2-2} and \ref{lem esp2-3}.
\end{Theorem}
\begin{Theorem}\label{th rho}
Under assumptions \emph{\textbf{(H.1-5)}} and if $\kappa\geq 8$, the following convergences hold
\begin{equation*}
\lim_{n\rightarrow\infty} \ind{|\dG_n^*|>0}{{\bs{\rho}}}_n={{\bs{\rho}}}\indnex\qquad a.s.
\end{equation*}
and
\begin{equation*}
\lim_{n\rightarrow\infty}\ind{|\dG_n^*|>0}\frac{|\dT_{n-1}^*|}{n}(\wh{{\bs{\rho}}}_n-{{\bs{\rho}}}_n)=\bs{V}^{-1}\big(q_{01}(0),  2q_{01}(1), q_{01}(2)\big)^t\indnex\qquad a.s.
\end{equation*}
where $\bs{V}$ is a $3\times3$ matrix defined in Proposition~\ref{prop lim SUV} and the $q_{01}(r)$ are scalars defined in Lemmas~\ref{lem P0}, \ref{lem P1} and \ref{lem P2}.
\end{Theorem}
We now turn to  the asymptotic normality for all our estimators
$\wh{\bs{\theta}}_n$, $\wh{\bs{\sigma}}_n$ and $\wh{\bs{\rho}}_n$ given the non-extinction of the underlying Galton-Watson process.  For this, using the fact that $\dP(\overline{\cE})\neq 0$ thanks to the super-criticality  assumption \textbf{(H.4)}, we define the probability $ \dP_{\overline{\cE}}$ on $(\Omega,\cA)$ by
$$ \dP_{\overline{\cE}}(A) =  \frac{\dP(A \cap \overline{\cE})}{ \dP(\overline{\cE})}\qquad \text{ for all } A \in {\cal A}.$$
\begin{Theorem}\label{thmCLT}
Under assumptions \emph{\textbf{(H.1-5)}} and if $\kappa\geq 4$, the following central limit theorem holds
\begin{equation}\label{CLTtheta}
{|\dT^*_{n-1}|}^{1/2} (\widehat{\bs{\theta}}_{n}-\bs{\theta})
\liml
\cN(0,\bs{S}^{-1}\bs{\Gamma} \bs{S}^{-1})\quad\text{on }\ 
({\overline{\cE}}, \dP_{\overline{\cE}})
\end{equation}with  $\bs{S}$ defined in Proposition \ref{prop lim SUV} and $\bs{\Gamma}$ in  Lemma~\ref{lem lim Gamma}.
If moreover $\kappa\geq 8$, 
\begin{equation}
\label{CLTsigma}
{|\dT^*_{n-1}|}^{1/2}  (\wh{\bs{\sigma}}_n-\bs{\sigma})
\liml
\cN\Bigl(0,\bs{U}^{-1}\bs{\Gamma^{\sigma}} \bs{U}^{-1} \Bigr)\quad\text{on }\ 
({\overline{\cE}}, \dP_{\overline{\cE}}),
\end{equation}and 
\begin{equation}
\label{CLTrho}
{|\dT^*_{n-1}|}^{1/2}  (\wh{\bs{\rho}}_n-\bs{\rho})
\liml
\cN(0,\bs{V}^{-1}\bs{\Gamma^{\rho}}\bs{V}^{-1})\quad\text{on }\ 
({\overline{\cE}}, \dP_{\overline{\cE}}),
\end{equation}with $\bs{U}$ and $\bs{V}$ defined in Proposition \ref{prop lim SUV} and $\bs{\Gamma^\sigma}$ and  $\bs{\Gamma^{\rho}}$ defined in Eq.~(\ref{Gammasigma}) and (\ref{Gammarho}).
\end{Theorem}
The proofs of these theorems are detailed in the next sections.

\section{Bifurcating Markov chains and consistency}
\label{section BMC}
In order to investigate the convergence of our estimators, we need laws of large numbers for quantities such as  $(\delta_{2k+i}X^q_kX^r_{2k}X_{2k+1}^s)_{k\in\dT}$. To obtain them, we use the bifurcating Markov chain framework introduced by J. Guyon in \cite{Guy07} and adapted to Galton-Watson trees by J.-F. Delmas and L. Marsalle in \cite{DM08}. We first recall the general framework, then prove the ergodicity of the induced Markov chain and finally derive strong laws of large numbers. We conclude this section by establishing the strong consistency of our estimators. Note that we cannot directly use the results in \cite{DM08} because our noise sequences do not have moments of all order. Therefore, our first step is to provide a general result for bifurcating Markov chains on GW trees with only a finite number of moments.

\subsection{Bifurcating Markov chain}
\label{subsection BMC}
Let $\mathbb{B}$ be the Borel $\sigma$-field of $\dR$, and $\mathbb{B}^p$ be the Borel $\sigma$-field of $\dR^p$, for $p>1$. We add a cemetery point $\partial$ to $\mathbb{R}$, denote by $\overline{\mathbb{R}}$ the set $\mathbb{R} \cup \{\partial \}$, and by $\overline{\mathbb{B}}$ the $\sigma$-field generated by $\mathbb{B}$ and $\{\partial\}$. This cemetery point models the state of a non-observed cell. We recall the following definitions from \cite{DM08}.

\begin{defi}
We call $\dT^*$-transition probability any mapping $P$ from $\overline{\dR} \times \overline{\dB}^2$ onto $[0,1]$ such that
\begin{itemize}
\item $P(\cdot,A)$ is measurable for all $A$ in $\overline{\dB}^2$,
\item $P(x,\cdot)$ is probability measure on $(\overline{\dR}^2,\overline{\dB}^2)$ for all $x$ in $\overline{\dR}$,
\item $P(\partial, \{(\partial, \partial)\}) = 1$.
\end{itemize}
\end{defi}
For any measurable function $f$ from $\overline{\dR}^3$ onto $\dR$, one defines the measurable function $Pf$ from $\overline{\dR}$ onto $\dR$ by
\begin{equation*}
Pf(x)=\int f(x,y,z)P(x,dy,dz),
\end{equation*}
provided the integral is well defined. Let $\nu$ be a probability measure on ${\dR}$. In the sequel, $\nu$ will denote the distribution of $X_1$.
\begin{defi}
We say that $(Z_n)_{n\in\dT}$ is a bifurcating Markov chain with initial distribution $\nu$ and $\dT^*$-transition probability $P$, a $P$-BMC in short, if
\begin{itemize}
\item $Z_1$ has distribution $\nu$,
\item for all $n$ in $\dN$, and for all family of measurable bounded functions $(f_k)_{k\in \dG_n}$ on $\dR^2$, 
\begin{equation*}
\dE\left[\prod_{k\in\dG_n}f_k(Z_{2k},Z_{2k+1}) \  \Big|\  \sigma(Z_j, j \in \dT_n)\right]=\prod_{k\in\dG_n}Pf_k(Z_k).
\end{equation*}
\end{itemize}
\end{defi} 
As explained in \cite{Guy07}, this means that given the first $n$ generations $\dT_n$, one builds generation $\dG_{n+1}$ by drawing $2^n$ independent couples $(Z_{2k},Z_{2k+1})$ according to $P(Z_k,\cdot)$, $k\in\dG_n$. This also means that any couple $(Z_{2k},Z_{2k+1})$ depends on past generations only through its mother $Z_k$.
The assumption $P(\partial, \{(\partial, \partial)\}) = 1$ means that $\partial$ is an absorbing state, and this hypothesis corresponds to the fact that a cell that is not observed cannot give birth to an observed one. We also assume that $P(x, \dR^2)$, $P(x, \dR \times \{\partial \})$ and $P(x, \{\partial\} \times \dR)$ do not depend on $x \in \dR$. The $P$-BMC is thus said to be spatially homogeneous. Such a spatially homogeneous $P$-BMC with an absorbing cemetery state is called a bifurcating Markov chain on a Galton-Watson tree, see \cite{DM08} for details.

Now let us turn back to our observed R-BAR process. In order to use the framework of $P$-BMC's, we define the auxilliary process $(X^*_n)_{n\in\dT}$ by 
\begin{equation}\label{Xstar}
X^*_n = X_n\ind{\delta_n=1} + \partial\ind{\delta_n=0},
\end{equation}
which means that $X_n^*=X_n$ if cell $n$ is observed, $X_n^*=\partial$ the cemetery state otherwise. It is clear from assumptions \textbf{(H.1)} and \textbf{(H.3)} that the process $(X^*_n)_{n\in\dT}$ is a P-BMC on a GW tree with $\dT^*$-transition probability given for all $x$ in $\dR$ and all measurable non-negative functions $f$ on $\overline{\dR}^3$ by
\begin{eqnarray}\label{def P}
Pf(x)&=&p_{01}\mathbb{E}\left[f\big(x,(b+\eta_2)x+a+\veps_2 ,(d+\eta_3)x+c+\veps_3\big)\right]
 + p_0\mathbb{E}\left[f\big(x,(b+\eta_2)x+a+\veps_2 ,\partial\big)\right]\nonumber \\
&&+ p_1\mathbb{E}\left[f\big(x,\partial ,(d+\eta_3)x+c+\veps_3\big)\right] + (1-p_{01}-p_0-p_1)f(x,\partial, \partial),
\end{eqnarray}
if $x\neq\partial$ and $Pf(\partial)=f(\partial,\partial, \partial)$.
As explained in \cite{Guy07}, the asymptotic behavior of the $P$-BMC is driven by that of the \emph{induced} Markov chain $(Y_n)$ defined on $\dR$ as follows. 
\begin{itemize}
\item For all $n\geq 1$, define  the sequence $(A_n,B_n)_{n\geq 1}$ to be i.i.d. random variables with the same distribution as $(a_{2+\zeta}+\veps_{2+\zeta}, b_{2+\zeta}+\eta_{2+\zeta})$, where $\zeta$ is a  Bernoulli random variable with mean $(p_{01}+p_1)/m$ independent from $(\veps_2,\eta_2,\veps_3,\eta_3)$.
\item Then, set $Y_0=X_1^*=X_1$ and define $Y_{n+1}$  recursively  by 
\begin{equation}\label{AR-Y}
Y_{n+1}=A_{n+1}+B_{n+1}Y_{n}.
\end{equation}
\end{itemize}
The sequence $(Y_ n)_{n\in\dN}$ is clearly an $\dR$-valued Markov chain with transition kernel given for all $x$ in $\dR$ and $A$ in $\dB$ by
\begin{equation}
\label{defQ}
Q(x,A) = \frac{P_0(x,A) + P_1(x,A)}{m}, 
\end{equation}
with $P_i(x,A)  =  (p_{01}+p_i)\mathbb{E}\left[\mathbbm{1}_{A}\big((b_{2+i}+\eta_{2+i})x+a_{2+i}+\veps_{2+i}\big)\right]$.
Note that the $P_0$ and $P_1$ are sub-probability kernels on ($\dR , \dB$), whereas $Q$ is a proper probability kernel on ($\dR , \dB$).
\subsection{Ergodicity of the induced Markov chain}
\label{subsection ergo}
We now turn to the investigation of ergodicity for the induced Markov chain $(Y_n)_{n\in\dN}$. We start with some preliminary results on the random variables $A_1$ and $B_1$.

\begin{Lemma}\label{lem AB}
Under assumptions \emph{\textbf{(H.{2})}} and \emph{\textbf{(H.5)}}, the random variables $A_1$ and $B_1$ have moments of all order up to $4\gamma$. In addition, $\mathbb{E}[\log|B_1|]<0$ and for all $0<s\leq 4\kappa$, the inequality $\mathbb{E}[|B_1|^s]<1$ holds.
\end{Lemma}

\noindent\textbf{Proof } First, for all $0\leq s\leq 4\gamma$, the following equalities clearly hold
\begin{equation*}
\mathbb{E}[|A_1|^s]=\frac{p_0+p_{01}}{m}\mathbb{E}[|a+\veps_2|^s]+\frac{p_1+p_{01}}{m}\mathbb{E}[|c+\veps_3|^s],\quad
\mathbb{E}[|B_1|^s]=\frac{p_0+p_{01}}{m}\mathbb{E}[|b+\eta_2|^s]+\frac{p_1+p_{01}}{m}\mathbb{E}[|d+\eta_3|^s].
\end{equation*}
Hence, under assumption \textbf{(H.2)}, it is clear that $\mathbb{E}[|A_1|^s]$ and $\mathbb{E}[|B_1|^s]$ are finite. {Next, notice that the function $s\longmapsto \mathbb{E}[|B_1|^s]$ is convex, that $\mathbb{E}[|B_1|^0]=1$ and $\mathbb{E}[|B_1|^{4\kappa}]<1$ by assumption \textbf{(H.5)}. This implies that $\mathbb{E}[|B_1|^s]<1$ for all $0<s\leq 4\kappa$. Last, consider $\mathbb{E}[|\log |B_1||]$: if it is finite, $\mathbb{E}[\log |B_1|]$ is the right-derivative at $0$ of $s\longmapsto \mathbb{E}[|B_1|^s]$, and convexity arguments with assumption \textbf{(H.5)} yield that  $\mathbb{E}[\log|B_1|]<0$ ; if it is infinite, the moment assumption on $B_1$ gives that necessarily $\mathbb{E}[(\log|B_1|)^+] < \infty$ and $\mathbb{E}[(\log|B_1|)^-] = \infty$, so that finally $\mathbb{E}[\log|B_1|]= -\infty <0$, as expected.}
\hspace{\stretch{1}}$ \Box$\\
 
\noindent The next result states the existence of an invariant distribution for the Markov chain $(Y_n)_{n\in\dN}$.  It is well known as $(Y_n)$ is a real-valued auto-regressive process with random i.i.d. coefficients satisfying Lemma~\ref{lem AB}, see e.g. \cite{Bra86,S05}.

\begin{Lemma}\label{lem mu existe}
Under assumptions \emph{\textbf{(H.{2})}} and \emph{\textbf{(H.5)}}, there exists a probability distribution $\mu$ on $(\dR,\dB)$ {(which is the distribution of the convergent series $Y_{\infty}=\sum_{\ell=1}^{\infty}B_1B_2\cdots B_{\ell-1}A_{\ell}$) such that} 
for all continuous bounded functions $f$ on $\dR$ and all $x$ in $\dR$, the following equality holds
\begin{equation*}
\dE_x[f(Y_n)]\xrightarrow [n\rightarrow\infty]{}\int fd\mu=\langle\mu,f\rangle.
\end{equation*}
\end{Lemma}
We investigate the moments of the invariant distribution $\mu$ to extend the above result to polynomial functions. For all $s\geq 1$, set $\|X\|_s=(\dE[|X|^s])^{1/s}$.

\begin{Lemma}\label{lem mu moment}
Under assumptions \emph{\textbf{(H.{2})}} and \emph{\textbf{(H.5)}}, $\mu$ has moments of all order up to $4\kappa$. In addition, for all $1\leq s\leq 4\kappa$, all $x\in \mathbb{R}$ and all $n\in \dN$, $(\dE_x[|Y_n|^s])^{1/s}\leq |x|+\|A_1\|_s/(1-\|B_1\|_s)<\infty$.
\end{Lemma}

\noindent\textbf{Proof } Set $1\leq s\leq 4\kappa$. As the sequence $(A_n,B_n)$ is i.i.d., the following inequality holds
\begin{eqnarray*}
\mathbb{E}[|Y_{\infty}|^s]^{1/s}&=&\dE\big[|\sum_{\ell=1}^{\infty}B_1\cdots B_{\ell-1}A_{\ell}|^s\big]^{1/s}\ \leq\ \sum_{\ell=1}^{\infty}\|B_1\|_s^{\ell-1}\|A_1\|_s.
\end{eqnarray*}
Since $\dE[|B_1|^s]<1$  and $\dE[|A_1|^s]<\infty$ thanks to Lemma~\ref{lem AB}, the series converges. Now let us turn to $Y_n$. The recursive equation~(\ref{AR-Y}) yields
\begin{equation*}
Y_n=Y_0B_1\cdots B_n+\sum_{\ell=1}^nB_n\cdots B_{\ell+1}A_{\ell},
\end{equation*}
with the usual convention that an empty product equals 1. As the sequence $(A_n,B_n)$ is i.i.d., $Y_n$ also has the same distribution (under $\dP_x$) as 
\begin{equation}\label{Yn}
xB_1\cdots B_n+\sum_{\ell=1}^nB_1\cdots B_{\ell-1}A_{\ell},
\end{equation}
so that, for $1\leq s\leq 4\kappa$, the following inequality holds
\begin{equation*}
\mathbb{E}_x[|Y_n|^s]^{1/s}\leq|x|\|B_1\|_s^n+\sum_{\ell=1}^n\|B_1\|_s^{\ell-1}\|A_1\|_s\leq|x|+\frac{\|A_1\|_s}{1-\|B_1\|_s},
\end{equation*}
hence the result.
\hspace{\stretch{1}}$ \Box$\\

\noindent The next result is a direct consequence of Lemma~\ref{lem mu moment}.

\begin{Corollary}\label{cor pol}
Under assumptions \emph{\textbf{(H.{2})}} and \emph{\textbf{(H.5)}}, all polynomial functions $f$ of degree less than $4\kappa$ are in $L_1(\mu)$: $\langle\mu,|f|\rangle=\dE[|f(Y_{\infty})|]<\infty$.
\end{Corollary}
We state a technical domination result that will be useful in the next section.

\begin{Lemma}\label{Q domine}
Under assumptions \emph{\textbf{(H.{2})}} and \emph{\textbf{(H.5)}}, for all polynomial functions $f$ of degree less than $2q$ with $q\leq 2\kappa$, there exists a nonnegative polynomial function $g$ of degree less than $2q$ such that for all $n\in\dN$ and all $x\in\dR$,
\begin{equation*}
\Big|\dE_x[f(Y_n)]\Big|\leq g(x).
\end{equation*}
\end{Lemma}

\noindent\textbf{Proof } It is sufficient to prove the result for $f(x)=x^p$ with $p \le 2q$. For $p \ge 1$, Lemma~\ref{lem mu moment} yields
\begin{eqnarray*}
\Big|\dE_x[Y_n^p]\Big| &\leq& \Big(|x|+\frac{\|A_1\|_p}{1-\|B_1\|_p}\Big)^p \le 2^{p-1}\Big(|x|^p+\frac{\|A_1\|_p^p}{(1-\|B_1\|_p)^p}\Big).
\end{eqnarray*}
If $p$ is even, we set $g(x)=2^{p-1}\Big(x^p+\|A_1\|_p^p/(1-\|B_1\|_p)^p\Big)$, and if $p$ is odd, we set $g(x)=2^{p-1}\Big(x^{p+1}+1+\|A_1\|_p^p/(1-\|B_1\|_p)^p\Big)$, as for all $x\in\dR$, $|x|^p\leq x^{p+1}+1$. Notice that if $p$ is odd and $p\leq 2q$, one also has $p+1\leq 2q$, hence the result.
\hspace{\stretch{1}}$ \Box$\\

\noindent Finally, we prove the geometric ergodicity of the induced Markov chain for polynomial functions.

\begin{Lemma}\label{Y ergodic}
Under assumptions \emph{\textbf{(H.{1-2})}} and \emph{\textbf{(H.5)}}, for all polynomial functions $f$ of degree less than $2q$ with $q\leq 2\kappa$, there exists a nonnegative polynomial function $g$ of degree less than $2q$ {and a positive constant $c$} such that for all $n\in\dN$ and all $x\in\dR$, the following inequalities hold
\begin{equation*}
\Big|\dE_x[f(Y_n)]-\langle\mu,f\rangle\Big|\leq g(x)\|B_1\|_{4\kappa}^n,\quad
{\Big|\dE_{\nu}[f(Y_n)]-\langle\mu,f\rangle\Big|\leq c \|B_1\|_{4\kappa}^n.}
\end{equation*}
\end{Lemma}

\noindent\textbf{Proof } Without loss of generality, it is sufficient to prove the result for polynomials $f$ of the form $x^p$ with $1\leq p\leq 2q$. H\"older's inequality yields
\begin{eqnarray*}
\Big|\dE_x[f(Y_n)]-\langle\mu,f\rangle\Big|&=&\Big|\dE_x[Y_n^p-Y_{\infty}^p]\Big|\ =\ \Big|\dE_x[(Y_n-Y_{\infty})\sum_{s=0}^{p-1}Y_n^sY_{\infty}^{p-1-s}]\Big|\\
&\leq&\Big(\dE_x[|Y_n-Y_{\infty}|^p]\Big)^{\frac{1}{p}}\sum_{s=0}^{p-1}\Big(\dE_x[|Y_n^{s}Y_{\infty}^{p-1-s}|^{\frac{p}{p-1}}]\Big)^{\frac{p-1}{p}}.
\end{eqnarray*}
 We are going to study the two terms above separately. For the first term, Eq.~(\ref{Yn}) and the definition of $Y_{\infty}$ yield
\begin{eqnarray*}
\Big(\dE_x[|Y_n-Y_{\infty}|^p]\Big)^{1/p}&=&\Big(\dE[|xB_1\cdots B_n-\sum_{\ell=n+1}^{\infty}B_1\cdots B_{\ell-1}A_{\ell}|^p]\Big)^{1/p}\\
&\leq&|x|\|B_1\|_p^n+\|A_1\|_p\frac{\|B_1\|_p^n}{1-\|B_1\|_p}\\
&\leq&\Big(|x|+\frac{\|A_1\|_p}{1-\|B_1\|_p}\Big)\|B_1\|_{4\kappa}^n,
\end{eqnarray*}
by Lemma~\ref{lem AB} as $p\leq 4\kappa$ by assumption. We now turn to the second term. H\"older's inequality with $\alpha=(p-1)/{s}$ and $\beta=(p-1)/(p-1-s)$ yields
\begin{eqnarray*}
\Big(\dE_x[|Y_n^{s}Y_{\infty}^{p-1-s}|^{\frac{p}{p-1}}]\Big)
&\leq&\Big(\dE_x[|Y_n|^{p}]\Big)^{s/(p-1)}\Big(\dE[|Y_{\infty}|^{p}]\Big)^{(p-1-s)/(p-1)} \\
&\leq&\Big(|x|+\frac{\|A_1\|_p}{1-\|B_1\|_p}\Big)^{s}\|Y_{\infty}\|_p^{p-1-s},
\end{eqnarray*}
this last upper bound coming from Lemma~\ref{lem mu moment}. Finally, one obtains
\begin{eqnarray*}
\Big|\dE_x[Y_n^p-Y_{\infty}^p]\Big|
&\leq& \|B_1\|_{4\kappa}^n \sum_{s=0}^{p-1} \Big(|x|+\frac{\|A_1\|_p}{1-\|B_1\|_p}\Big)^{s+1}\|Y_{\infty}\|_p^{p-1-s}\ \leq\ \|B_1\|_{4\kappa}^ng(x),
\end{eqnarray*}
where $g$ is a polynomial function of degree at most $2q$ by a similar argument as in the previous proof. {Integrating this upper bound with respect to the initial law $\nu$ and using \textbf{(H.1)} gives the second result.}
\hspace{\stretch{1}}$ \Box$

\subsection{Laws of large numbers for the $P$-BMC}
\label{subsection LLN}
We now want to prove laws of large numbers for a family of functionals of the $P$-BMC $(X_n^*)$. 
We are interested in polynomial functions on $\overline{\dR}$ and $\overline{\dR}^3$ multiplied by indicators. Precisely, for all $q \ge 1$, let $F_{q}$ and $G_q$ be the vector spaces generated by the following class of functions from $\overline{\dR}^3$ to $\overline{\dR}$ and from $\overline{\dR}$ to $\dR$ respectively,
$$\{x^\alpha  y^\beta \inds{\dR}(y),~x^\alpha  z^\tau \inds{\dR}(z),~x^\alpha y^\beta z^\tau \inds{\dR^2}(y,z),\quad 0\leq \alpha+\beta+\tau\leq q\},$$
$$\{x^\alpha \inds{\dR}(x),\quad 0\leq \alpha \leq q\},$$
{where $\alpha$, $\beta$, $\tau$ are integers.}
We first establish some technical results needed in the main proof. 

\begin{Lemma}\label{lem Pf}
Let $f \in F_q$ and $h \in G_q$. Under assumption \emph{\textbf{(H.{2})}}, 
\begin{itemize}
 \item[(i)] if $q\leq 4\gamma$, $f \in L^1(P)$ and $Pf \in G_q$,
 \item[(ii)] if $q\leq 4\gamma$, $h \in L^1(P_0,P_1,Q)$ and $P_0h, P_1h$ and $Qh \in G_q$,
 \item[(iii)] if $q\leq 2\gamma$, $h \otimes h \in L^1(P)$ and $P(h \otimes h) \in G_{2q}$.
\end{itemize}
\end{Lemma}

\noindent\textbf{Proof }
Take $q \le 4\gamma$ and $p \le 2\gamma$ and remark that $ Pf(\partial)=0$ for any $f \in F_q$, so that $Pf(x)=Pf(x) \inds{\dR}(x)$ for all $x \in \overline{\dR}$. Next, take $f_1(x,y,z)=x^\alpha y^\beta z^\tau\inds{\dR^2}(y,z)$ and $f_2=x^\delta y^\epsilon \inds{\dR}(y)$ in $F_q$, $h(x)=x^{\rho} \inds{\dR}(x)$ in $G_q$, and  $l(y,z)=y^i \otimes z^j \inds{\dR^2}(y,z)$ with $i,j \le p$. Eq.~(\ref{def P}) yields, for $i \in \{0,1\}$,
\begin{eqnarray*}
P|f_1|(x) & = & p_{01} |x|^{\alpha} \dE\big[\big|(b+\eta_2)x + a + \veps_2 \big|^{\beta} \big|(d+\eta_3)x+c+\veps_3 \big|^{\tau}\big],\\
P|f_2|(x) & = & (p_{01}+p_0) |x|^{\delta} \dE\big[\big|(b+\eta_2)x + a + \veps_2 \big|^{\epsilon}\big],\\
P_i|h|(x) & = & (p_{01} + p_i)  \dE\big[\big|(b_{2+i}+\eta_{2+i})x + a_{2+i} + \veps_{2+i} \big|^{\rho} \big],\\
P|l|(x) & = & p_{01}  \dE\big[\big|(b+\eta_2)x + a + \veps_2 \big|^i \big|(d+\eta_3)x+c+\veps_3 \big|^j\big]. 
\end{eqnarray*}
Assumption {\textbf{(H.2)}} entails that the $4\gamma$-th moments of $\big((b+\eta_2)x + a + \veps_2\big)$ and $\big((d+\eta_3)x+c+\veps_3\big)$ are finite, which gives the integrability results, since $\beta + \tau \le q$, $\epsilon \le q$, $\rho \le q$ and $i+j \le 2p$. The integrability results are thus proved. It is then obvious that $Pl(x)$ is computed the same way as $Pf_1(x)$, and $P_ih(x)$ the same way as $Pf_2(x)$. But
\begin{eqnarray*}
Pf_1(x) &=& p_{01}\sum_{r=0}^{\beta} \sum_{s=0}^\tau C_\beta^r C_\tau^s \dE\big[(b+\eta_2)^r(a + \veps_2)^{\beta-r} (d+\eta_3)^s  (c+\veps_3)^{\tau-s}\big] x^{r+s+\alpha},\\
Pf_2(x) &=&  (p_{01}+p_0)\sum_{r=0}^{\epsilon} C_\epsilon^r \dE\big[(b+\eta_2)^r(a + \veps_2)^{\epsilon-r}\big] x^{r+\delta}.
\end{eqnarray*}
so that $Pf_1$, $Pf_2$ and $P_ih$ are in $G_q$, since $\alpha +\beta+\tau \le q$, $\delta+\epsilon \le q$ and $\rho \le q$. As for $Pl$, it belongs to $G_{2p}$, since $i+j \le 2p$.
\hspace{\stretch{1}}$ \Box$\\

\noindent We are now ready to prove the main result of this section.

\begin{Theorem}\label{th LGN}
Under assumptions \emph{\textbf{(H.1-5)}},  for all function $f \in F_{\kappa}$,  the following law of large numbers holds
\begin{equation*}
\lim_{n\rightarrow\infty}\frac{1}{m^n} \sum_{k \in \dG_n^*} f(\bX_k, \bX_{2k}, \bX_{2k+1}) =\langle \mu, Pf\rangle W \quad \text{a.s.}
\end{equation*}
\end{Theorem}

\noindent\textbf{Proof: } This result is similar to Theorem~11 of \cite{Guy07} and Theorem~3.1 of \cite{DM08}. The proof follows essentially the same lines and is thus shortened here, the main difference being that the class of functions $F_{\kappa}$ does not satisfy assumptions (i)-(vi) from \cite{Guy07, DM08} mainly because $F_{\kappa}$ is not stable by multiplication and $(\veps_2,\eta_2, \veps_3, \eta_3)$ do not have moments of all order.

 For all $f$ in $F_{\kappa}$, $Pf$ is well-defined from $\overline{\dR}$ onto $\dR$ thanks to Lemma~\ref{lem Pf} as $\kappa\leq\gamma$. As $Pf(\partial)=0$, by a slight abuse of notation we will also denote $Pf$ its restriction to $\dR$. Thus, $Pf$ is $\mu$-integrable by Lemma \ref{lem mu moment}. One has
\begin{eqnarray*}
{m^{-n} \sum_{k \in \dG_n^*} f(\bX_k, \bX_{2k}, \bX_{2k+1}) - \langle \mu, Pf\rangle W}
& = & \frac{1}{m^n}\sum_{k \in \dG_n^*} \big(f(\bX_k, \bX_{2k}, \bX_{2k+1}) - \langle \mu, Pf\rangle \big) + \langle \mu, Pf\rangle \big(\frac{|\dG_n^*|}{m^n} -W\big).
\end{eqnarray*} By Eq.~(\ref{eq cvGW}) the second term converges to 0 a.s. as $n$ tends to infinity. In order to prove the a.s. convergence of the first term, as in \cite{Guy07, DM08}, it is sufficient to prove that
\begin{equation}\label{eq:limps}
\sum_{n \ge 0} m^{-2n} \dE\big[\big(\sum_{k \in \dG_n^*} g(\bX_k, \bX_{2k}, \bX_{2k+1})\big)^2\big] < \infty,
\end{equation}
with $g=f-\langle \mu, Pf\rangle\in F_{\kappa}$. Thanks to Lemma~\ref{lem Pf},  $Pg \in G_{\kappa}$, and as $g^2 \in F_{2\kappa}$, one also has $Pg^2\in G_{2\kappa}$. 
The expectation inside the sum decomposes as
\begin{eqnarray*}
{\dE\big[\big(\sum_{k \in \dG_n^*} g(\bX_k, \bX_{2k}, \bX_{2k+1})\big)^2\big] }
&=& {\dE\big[\big(\sum_{k \in \dG_n^*} Pg(\bX_k)\big)^2\big]} + {\dE\big[\sum_{k \in \dG_n^*} \big(Pg^2 - (Pg)^2\big)(\bX_k) \big]}\ = \ C_n+D_n.
\end{eqnarray*}
We study the two terms $C_n$ and $D_n$ separately. Let us first prove that $\sum_{n \ge 0}  m^{-2n} D_n < \infty$. We can rewrite $D_n = \dE\big[\sum_{k \in \dG_n^*} h(\bX_k)\big]$ with $h = Pg^2-(Pg)^2$. As seen above, $h\in G_{2\kappa}$ and therefore $h$ is $\mu$-integrable thanks to Lemma~\ref{lem mu moment}. To investigate the limit of $\sum m^{-2n}D_n$, we prove that $m^{-n}D_n$ has a finite limit. More precisely, the following inequality holds
\begin{eqnarray*}
{||m^{-n}\sum_{k \in \dG_n^*} h(\bX_k) - \langle \mu, h \rangle W||_{2}}
& =& ||m^{-n}\sum_{k \in \dG_n^*} \big(h(\bX_k) - \langle \mu , h \rangle \big) +  \langle \mu , h \rangle (m^{-n}|\dG_n^*|-W) ||_{2}\\
&\le& ||m^{-n}\sum_{k \in \dG_n^*} \big(h(\bX_k) - \langle \mu , h \rangle \big)||_{2} + |\langle \mu , h \rangle| ~~|| m^{-n}|\dG_n^*|-W||_{2}.
\end{eqnarray*}
The second term converges to zero. For the first term, again let  $l := h - \langle \mu , h \rangle \in G_{2\kappa}$ and $\langle \mu , l\rangle =0$, and by \cite[equation (15) p 2504]{DM08}, the following equality holds
\begin{eqnarray}\label{eq:cle}
||m^{-n}\sum_{k \in \dG_n^*} \big(h(\bX_k) - \langle \mu , h \rangle \big)||_{2}^2
& = & m^{-n} \dE_{\nu}[l^2(Y_n)] + 2m^{-2} \sum_{\ell=0}^{n-1} m^{-\ell} \langle \nu, Q^{\ell} P (Q^{n-\ell-1}l \otimes Q^{n-\ell-1}l) \rangle. 
\end{eqnarray}
Concerning the first term in Eq. (\ref{eq:cle}), as $l^2\in G_{4\kappa}$, by Lemma \ref{Y ergodic} one obtains $\lim_{n\rightarrow\infty} \dE_{\nu}[l^2(Y_n)] = \langle \mu , l^2 \rangle$ and $m^{-n} \dE[l^2(Y_n)] $ converges to 0 a.s.
Concerning the second term in Eq. (\ref{eq:cle}), Lemma~\ref{Y ergodic} yields 
$\lim_{n\rightarrow\infty} Q^{n-\ell-1}l(x) = \lim_{n\rightarrow\infty} \dE_x[l(Y_{n-\ell-1}]=\langle \mu,l\rangle=0$
 and by Lemma~\ref{Q domine}, $Q^{n-r-1}l$ is dominated by some $\phi \in G_{2\kappa}$. Moreover, using Lemma~\ref{lem Pf},
 $\phi \otimes \phi$ belongs to $F_{4\kappa}$, it is $P$-integrable and $P(\phi \otimes \phi)$ belongs to $G_{4\kappa}$.
 By Lemma~\ref {Q domine}, $Q^{\ell} P(\phi \otimes \phi)$ is dominated by  some $\psi \in G_{4\kappa}$, which is $\nu$-integrable by assumption \textbf{(H.2)}. Lebesgue dominated convergence theorem yields 
 $$ \lim_{n\rightarrow\infty} \langle \nu, Q^{\ell} P (Q^{n-\ell-1}l \otimes Q^{n-\ell-1}l) \rangle = 0,$$
and  $|\langle \nu, Q^{\ell} P (Q^{n-\ell-1}l \otimes Q^{n-\ell-1}l) \rangle| \le \langle \nu , \psi \rangle$. 
This upper bound allows us  to deal with the limit of  the second term of Eq. (\ref{eq:cle}). Under assumption {\textbf{(H.4)}}, $\sum_{\ell=0}^{n-1} m^{-\ell} $ converges and for $\eps >0$, there exists $\ell_{\epsilon}$ such that 
$ \sum_{\ell=\ell_ \epsilon}^{n-1} m^{-\ell} \langle \nu,\psi\rangle \leq \epsilon$.  
Finally, for $n > \ell_{\epsilon}$, we obtain
\begin{eqnarray*}
{\Big|\sum_{\ell=0}^{n-1} m^{-\ell} \langle \nu, Q^{\ell} P (Q^{n-\ell-1}l \otimes Q^{n-\ell-1}l) \rangle \Big| } 
& \le & \sum_{\ell=0}^{\ell_\epsilon-1} m^{-\ell} |\langle \nu, Q^{\ell} P (Q^{n-\ell-1}l \otimes Q^{n-\ell-1}l) \rangle| + \epsilon, 
\end{eqnarray*}
All the terms of the left sum converge to 0 with $n$, which finally proves the $L_2$-convergence of $m^{-n}\sum_{k \in \dG_n^*} h(\bX_k)$ to $\langle \mu, h\rangle W$. This implies the convergence of the expectation $ m^{-n} D_n$ to $\langle \mu, h\rangle\dE[W]$ (recall that $W$ is square-integrable). Therefore, one obtains $\sum_{n \ge 0}  m^{-2n} D_n < \infty$ because $m>1$.

Let us now prove that $\sum_{n \ge 0}  m^{-2n} C_n < \infty$. Recall that $g \in F_{\kappa}$, $\langle \mu, Pg \rangle = 0$ and following Eq.~(15) p. 2504 of \cite{DM08}, we obtain a new expression for $C_n$: 
\begin{eqnarray}
\frac{C_n} {m^{2n}}& = &  ||\frac{1}{m^n}\sum_{k \in \dG_n^*} Pg(X_k^*)||_{2}^2 \label{eq:Aq} \\
& = & \frac{1}{m^n} \dE_{\nu}[(Pg)^2(Y_n)] + \frac{2}{m^2} \sum_{\ell=0}^{n-1}\frac{\langle \nu, Q^{\ell} P (Q^{n-\ell-1}(Pg) \otimes Q^{n-\ell-1}(Pg)) \rangle}{m^{\ell}} . \nonumber
\end{eqnarray}
The proof of the convergence of the first term of Eq.~(\ref{eq:Aq}) is the same as the one of $\dE_{\nu}[l^2(Y_q)] $, and $\sum_{n \ge 0} m^{-n} \dE[(Pg)^2(Y_n)]$ converges. For the second term, setting $p=n-\ell-1$, we can rewrite
\begin{eqnarray*}
{\sum_{n \ge 0} \sum_{\ell=0}^{n-1} m^{-\ell} \langle \nu, Q^{\ell} P (Q^{n-\ell-1}(Pg) \otimes Q^{n-\ell-1}(Pg)) \rangle} 
& = & \sum_{\ell\ge 0} m^{-\ell} \langle \nu, Q^{\ell} P \Big(\sum_{p \ge 0} (Q^{p}(Pg) \otimes Q^{p}(Pg))\Big) \rangle. \label{pq'}
\end{eqnarray*}
By Lemma \ref{Y ergodic}, there exists $\varphi \in G_{\kappa+1}$, such that 
$|\dE_x[(Pg)(Y_{p})]|=|Q^{p}(Pg)(x)| \le \varphi(x)  \|B_1\|_{4\kappa}^{p}$
and therefore the following inequality holds
$$|\sum_{p \ge 0} (Q^{p}(Pg) \otimes Q^{p}(Pg))| \le (\varphi \otimes\varphi) \sum_{p \ge 0} \|B_1\|_{4\kappa}^{p}.$$
By assumption~{\textbf{(H.5)}}, the series converges and there only remains to study the asymptotic behavior of
$\sum_{\ell\ge 0} m^{-\ell} \langle \nu, Q^{\ell} P(\varphi \otimes \varphi) \rangle.$ For this, let us remark that
 $\langle \nu, Q^{\ell} P(\varphi \otimes \varphi) \rangle = \dE_{\nu}[P(\varphi \otimes \varphi)(Y_{\ell})]$ with  $P(\varphi \otimes \varphi) \in G_{2\kappa+2}$. By Lemma \ref{Y ergodic}, $\lim_{\ell \rightarrow\infty}\dE_{\nu}[P(\varphi \otimes \varphi)(Y_{\ell})]$ is finite and the series converges because $m > 1$.
We have thus proved that Eq.~(\ref{eq:limps}) holds, and hence the almost sure convergence of the series ${m^{-n}} \sum_{\ell \in \dG_n^*} f(\bX_{\ell}, \bX_{2\ell}, \bX_{2\ell+1})$ to $\langle \mu, Pf\rangle W$.
\hspace{\stretch{1}}$ \Box$

\subsection{Laws of large numbers for the R-BAR process}
\label{subsection LLNMC}
Let us now turn back to our R-BAR process and see how the law of large numbers given by Theorem~\ref{th LGN} applies to our process.

\begin{Proposition}\label{prop LGN}
Under assumptions \emph{\textbf{(H.1-5)}},  for all integers $0\leq q\leq {\kappa}$, and all $i\in\{0,1\}$, the following laws of large numbers hold
\begin{eqnarray*}
\lim_{n\rightarrow\infty}\ind{|\dG_n^*|>0}\frac{1}{|\dT_n^*|} \sum_{k \in \dT_n^*} \delta_{2k+i}X_k^q& =&\ell_i(q)\indnex\quad \text{a.s.}\\
\lim_{n\rightarrow\infty}\ind{|\dG_n^*|>0}\frac{1}{|\dT_n^*|} \sum_{k \in \dT_n^*} \delta_{2k}\delta_{2k+1}X_k^q &=&\ell_{01}(q)\indnex\quad \text{a.s.}
\end{eqnarray*}
with $\ell_i(q)=(p_{01}+p_i)\dE[Y_{\infty}^q]$ and $\ell_{01}(q)=p_{01}\dE[Y_{\infty}^q]$.
\end{Proposition}

\noindent\textbf{Proof } Set $q\leq \kappa$. We apply Theorem~\ref{th LGN} to the function $f_0(x,y,z)=x^q\inds{\dR}(y)$ if $i=0$ and $f_1(x,y,z)=x^q\inds{\dR}(z)$ if $i=1$ for the first limit, and $f_{01}(x,y,z)=x^q\inds{\dR^2}(y,z)$ for the second limit. The functions $f_0$, $f_{1}$ and $f_{01}$ clearly belong to $F_{\kappa}$, and moreover $Pf_i(x)=(p_{01}+p_i)x^q$, $Pf_{01}(x)=p_{01}x^q$. Finally, notice that $\langle \mu,x^q\rangle=\dE[Y_{\infty}^q]$. Theorem~\ref{th LGN} thus yields 
\begin{equation*}
\lim_{n\rightarrow\infty}\frac{1}{m^n} \sum_{k \in \dG_n^*} \delta_{2k+i}X_k^q =\ell_i(q)W,\quad
\lim_{n\rightarrow\infty}\frac{1}{m^n} \sum_{k \in \dG_n^*} \delta_{2k}\delta_{2k+1}X_k^q =\ell_{01}(q)W\quad \text{a.s.}
\end{equation*}
Now, for instance, the following decomposition holds
\begin{equation*}
\lim_{n\rightarrow\infty}\frac{1}{m^n} \sum_{k \in \dT_n^*} \delta_{2k+i}X_k^q=\sum_{\ell=0}^n\frac{1}{m^{n-\ell}}\Big(\frac{1}{m^{\ell}}\sum_{k\in\dG_{\ell}^*}\delta_{2k+i}X_k^q\Big).
\end{equation*}
The sum above converges to $\ell_i(q)Wm/(m-1)$ thanks to Lemma A.3 of \cite{BSG09} and we conclude 
using Eq.~(\ref{eq cvGW2}).
\hspace{\stretch{1}}$ \Box$

\begin{Proposition}\label{prop LGN2}
Under assumptions \emph{\textbf{(H.1-5)}}, for all integers $0\leq q\leq {\kappa}-1$, and all $i\in\{0,1\}$,  the following almost sure convergences hold
\begin{equation*}
\lim_{n\rightarrow\infty}\frac{\ind{|\dG_n^*|>0}}{|\dT_{n-1}^*|} \sum_{k \in \dT_{n-1}^*} \delta_{2k+i}X_k^qX_{2k+i} =(p_{01}+p_i)(a_{2+i}\dE[Y_{\infty}^q]+b_{2+i}\dE[Y_{\infty}^{q+1}])\indnex,
\end{equation*}
\begin{equation*}
\lim_{n\rightarrow\infty}\frac{\ind{|\dG_n^*|>0}}{|\dT_{n-1}^*|} \sum_{k \in \dT_{n-1}^*} \delta_{2k}\delta_{2k+1}X_k^qX_{2k+i} =p_{01}(a_{2+i}\dE[Y_{\infty}^q]+b_{2+i}\dE[Y_{\infty}^{q+1}])\indnex,
\end{equation*}
and if $\kappa\geq 2$, for all integers $0\leq q\leq {\kappa}-2$,  the following almost sure convergences holds
\begin{eqnarray*}
\lefteqn{\lim_{n\rightarrow\infty}\ind{|\dG_n^*|>0}\frac{1}{|\dT_{n-1}^*|} \sum_{k \in \dT_{n-1}^*} \delta_{2k+i}X_k^qX_{2k+i}^2}\\
& =&(p_{01}+p_i)\big((a_{2+i}^2+\sigma^2_{\veps})\dE[Y_{\infty}^q]+2(a_{2+i}b_{2+i}+\rho_{ii})\dE[Y_{\infty}^{q+1}]+(b_{2+i}^2+\sigma_{\eta}^2)\dE[Y_{\infty}^{q+2}]\big)\indnex,\\
\lefteqn{\lim_{n\rightarrow\infty}\ind{|\dG_n^*|>0}\frac{1}{|\dT_{n-1}^*|} \sum_{k \in \dT_{n-1}^*} \delta_{2k}\delta_{2k+1}X_k^qX_{2k}X_{2k+1}}\\
& =&p_{01}\big((ac+\rho_{\veps})\dE[Y_{\infty}^q]+(ad+bc+2\rho)\dE[Y_{\infty}^{q+1}]+(bd+\rho_{\eta})\dE[Y_{\infty}^{q+2}]\big)\indnex.
\end{eqnarray*}
\end{Proposition}

\noindent\textbf{Proof } The proof follows the same lines as that of Proposition~\ref{prop LGN}.
\hspace{\stretch{1}}$ \Box$\\

We end this section by stating how to compute the moments of the invariant law $\mu$.

\begin{Lemma}\label{lem LGN}
Under assumptions \emph{\textbf{(H.2)}} and \emph{\textbf{(H.5)}}, the first moments of $Y_{\infty}$ are
\begin{equation*}
\dE[Y_{\infty}]=\frac{\dE[A_1]}{1-\dE[B_1]},\qquad \dE[Y_{\infty}^2]=\frac{\dE[A_1^2]+2\dE[A_1B_1]\dE[Y_{\infty}]}{1-\dE[B_1^2]},
\end{equation*}
and more generally, the moments of $Y_{\infty}$ can be calculated recursively for all $1\leq q\leq 4\kappa$ thanks to the relation
$\dE[Y_{\infty}^q]=\sum_{s=0}^qC_q^s\dE[A_1^{q-s}B_1^{s}]\dE[Y_{\infty}^s]$.
\end{Lemma}

\noindent\textbf{Proof } As $Y_{\infty}$ is the stationary solution of equation $Y_n=A_n+B_nY_{n-1}$, $Y_{\infty}$ has the same law as $A_0+B_0Y_{\infty}$ where $(A_0,B_0)$ is a copy of $(A_1,B_1)$ independent from the sequence $(A_n,B_n)_{n\geq 1}$. Hence, we can write $\dE[Y_\infty]=\dE[A_0+B_0Y_{\infty}]=\dE[A_1]+\dE[B_1]\dE[Y_{\infty}]$. Similarly, one has 
\begin{eqnarray*}
\dE[Y_{\infty}^2]&=&\dE[(A_0+B_0Y_{\infty})^2]\ =\ \dE[A_1^2]+2\dE[A_1B_1]\dE[Y_{\infty}]+\dE[B_1^2]\dE[Y_{\infty}^2].
\end{eqnarray*}
The general formula is obtained in the same way by developing the relation $\dE[Y_{\infty}^q]=\dE[(A_0+B_0Y_{\infty})^q]$.
\hspace{\stretch{1}}$ \Box$\\

Note that one can easily compute the moments of $A_1$ and $B_1$ from their definition. In particular, the two following equalities hold
\begin{equation*}
\dE[Y_{\infty}]=\frac{am_0+cm_1}{1-bm_0-dm_1},
\end{equation*}
and 
\begin{equation*}
\dE[Y_{\infty}^2]=\frac{a^2m_0+c^2m_1+\sigma_{\veps}^2+2\big((ab+\rho_{00})m_0+(cd+\rho_{11})m_1\big)\dE[Y_{\infty}]}{1-(b^2m_0+d^2m_1+\sigma_{\eta}^2)},
\end{equation*}
with $m_0=(p_{01}+p_0)/m$ and $m_1=(p_{01}+p_1)/m$.
\subsection{Consistency of the estimators}
\label{subsection consistency}
We are now able to prove the consistency of our estimators. We start with the computation of the limits of the normalizing matrices $\bs{S}_n$, $\bs{U}_n$ and $\bs{V}_n$, which is a direct consequence of Proposition~\ref{prop LGN}
\begin{Proposition}\label{prop lim SUV}
Under assumptions \emph{\textbf{(H.1-5)}}, and if $\kappa\geq 2$, for $i\in\{0,1\}$, the following laws of large numbers hold
\begin{equation*}
\lim_{n\rightarrow\infty}\ind{|\dG_n^*|>0}\frac{\bs{S}^i_n}{|\dT_{n}^*|}=\bs{S}^i\indnex
=(p_{01}+p_i)\left(\begin{array}{cc}1&\dE[Y_{\infty}]\\\dE[Y_{\infty}]&\dE[Y_{\infty}^2]\end{array}\right)\indnex\quad \text{a.s.}
\end{equation*}and 
\begin{equation*}
\lim_{n\rightarrow\infty}\ind{|\dG_n^*|>0}\frac{\bs{S}_n}{|\dT_{n}^*|}=\bs{S}\indnex
=\left(\begin{array}{cc}
\bs{S}^0&0\\0&\bs{S}^1
\end{array}\right)\indnex\quad \text{a.s.}
\end{equation*}
If in addition $\kappa\geq 4$, the following convergences hold
\begin{eqnarray*}
{\lim_{n\rightarrow\infty}\ind{|\dG_n^*|>0}\frac{\bs{U}_n}{|\dT_{n}^*|}\ =\ \bs{U}\indnex}
&=&m\left(\begin{array}{rrrr}
1&2m_0\dE[Y_{\infty}]&2m_1\dE[Y_{\infty}]&\dE[Y_{\infty}^2]\\
2m_0\dE[Y_{\infty}]&4m_0\dE[Y_{\infty}^2]&0&2m_0\dE[Y_{\infty}^3]\\
2m_1\dE[Y_{\infty}]&0&4m_1\dE[Y_{\infty}^2]&2m_1\dE[Y_{\infty}^3]\\
\dE[Y_{\infty}^2]&2m_0\dE[Y_{\infty}^3]&2m_1\dE[Y_{\infty}^3]&\dE[Y_{\infty}^4]
\end{array}\right)\indnex\quad \text{a.s.}
\end{eqnarray*}
and
\begin{equation*}
\lim_{n\rightarrow\infty}\ind{|\dG_n^*|>0}\frac{\bs{V}_n}{|\dT_{n}^*|}=\bs{V}\indnex
=p_{01}\left(\begin{array}{rrr}
1&2\dE[Y_{\infty}]&\dE[Y_{\infty}^2]\\
2\dE[Y_{\infty}]&4\dE[Y_{\infty}^2]&2\dE[Y_{\infty}^3]\\
\dE[Y_{\infty}^2]&2\dE[Y_{\infty}^3]&\dE[Y_{\infty}^4]
\end{array}\right)\indnex\quad \text{a.s.}
\end{equation*}
Besides, the matrices $\bs{S}^i$, $\bs{U}$ and $\bs{V}$ are invertible. 
\end{Proposition}
We now turn to the consistency of our main estimators.\\

\noindent\textbf{Proof of Theorem~\ref{th consist}} As regards our main estimator $\bs{\wh{\theta}}_n$, a direct application of Proposition~\ref{prop LGN2} yields
\begin{equation*}
\lim_{n\rightarrow\infty}\frac{\ind{|\dG_n^*|>0}}{|\dT_{n-1}^*|}
\bs{S}_{n-1}\wh{\bs{\theta}}_n
=m\left( \begin{array}{c}
m_0(a+b\dE[Y_{\infty}])  \\
m_0(a\dE[Y_{\infty}]+b\dE[Y_{\infty}^2]) \\
m_1(c+d\dE[Y_{\infty}])  \\
m_1(c\dE[Y_{\infty}]+d\dE[Y_{\infty}^2]) 
\end{array}\right)\indnex
=\bs{S}\bs{\theta}\indnex
\ \text{a.s.}
\end{equation*}
and the result follows from Proposition~\ref{prop lim SUV} and the definition of $\bs{\wh{\theta}}_n$. 
The consistency of $\bs{\wh{\sigma}}_n$ and $\bs{\wh{\rho}}_n$ is a bit more complicated because their definition involves the $\wh{\eps}_k$. We give a detailed proof of the convergence of $|\dT_{n-1}^*|^{-1}\sum\wh{\eps}_{2k}^2$, the other terms in $\bs{U}_{n-1}\bs{\wh{\sigma}}_n$ and $\bs{V}_{n-1}\bs{\wh{\rho}}_n$ being treated similarly. For $k\in\dG_n$, one can develop
\begin{eqnarray*}
\wh{\eps}_{2k}^2&=&\delta_{2k}(X_{2k}-\wh{a}_n-\wh{b}_nX_k)^2\\
&=&\delta_{2k}(\wh{a}_n^2+2\wh{a}_n\wh{b}_nX_k+\wh{b}_n^2X_k^2-2\wh{a}_nX_{2k}-2\wh{b}_nX_kX_{2k}+X_{2k}^2).
\end{eqnarray*}
Hence, the following equality holds
\begin{eqnarray}
\sum_{k\in\dT_{n-1}^*}\wh{\eps}_{2k}^2
&=&\sum_{\ell=1}^{n-1}\wh{a}_{\ell}^2\sum_{k\in\dG_{\ell}}\delta_{2k}+2\sum_{\ell=1}^{n-1}\wh{a}_{\ell}\wh{b}_{\ell}\sum_{k\in\dG_{\ell}}\delta_{2k}X_k+\sum_{\ell=1}^{n-1}\wh{b}_{\ell}^2\sum_{k\in\dG_{\ell}}\delta_{2k}X_k^2\label{eq eps}\\
&&-2\sum_{\ell=1}^{n-1}\wh{a}_{\ell}\sum_{k\in\dG_{\ell}}\delta_{2k}X_{2k}-2\sum_{\ell=1}^{n-1}\wh{b}_{\ell}\sum_{k\in\dG_{\ell}}\delta_{2k}X_kX_{2k}+\sum_{k\in\dT_{n-1}^*}\delta_{2k}X_{2k}^2.\nonumber
\end{eqnarray}
The limit of the last term is given by Proposition~\ref{prop LGN2}. The first term decomposes as
\begin{equation*}
\frac{1}{m^{n-1}}\sum_{\ell=1}^{n-1}\wh{a}_{\ell}^2\sum_{k\in\dG_{\ell}}\delta_{2k}=\sum_{\ell=1}^{n-1}\wh{a}_{\ell}^2\frac{m^{\ell}}{m^{n-1}}\frac{1}{m^{\ell}}\sum_{k\in\dG_{\ell}}\delta_{2k}.
\end{equation*}
We apply Lemma A.3 of \cite{BSG09} to the sequence above. On the one hand,
\begin{equation*}
\lim_{\ell\rightarrow\infty}\wh{a}_{\ell}^2\frac{1}{m^{\ell}}\sum_{k\in\dG_{\ell}}\delta_{2k}=a^2(p_{01}+p_0)W\quad \text{a.s.}
\end{equation*}
 thanks to the previous result on the consistency of $\bs{\wh{\theta}}_n$ and Theorem~\ref{th LGN}. On the other hand, the series $\sum m^{-n}$ converges to $m/(m-1)$ under assumption \textbf{(H.4)}. Therefore, Lemma~A.3 of \cite{BSG09} yields
 \begin{equation*}
\lim_{n\rightarrow\infty}\frac{1}{m^{n-1}}\sum_{\ell=1}^{n-1}\wh{a}_{\ell}^2\sum_{k\in\dG_{\ell}}\delta_{2k}= \frac{m}{m-1}a^2(p_{01}+p_0)W\quad \text{a.s.}
 \end{equation*}
 and Eq.~(\ref{eq cvGW2}) finally yields 
 \begin{equation*}
\lim_{n\rightarrow\infty}\ind{|\dG_n^*|>0} \frac{1}{|\dT_{n-1}^*|}\sum_{\ell=1}^{n-1}\wh{a}_{\ell}^2\sum_{k\in\dG_{\ell}}\delta_{2k}= a^2(p_{01}+p_0)\indnex\quad \text{a.s.}
 \end{equation*}
 Note that the limit above is just the limit of $\wh{a}_{\ell}^2$ multiplied by the limit of $|\dT_{n-1}^*|^{-1}\sum\delta_{2k}$. The other terms in Eq.~(\ref{eq eps}) are dealt with similarly using the results of Proposition~\ref{prop LGN2}. Finally, one obtains the almost sure convergences
\begin{eqnarray*}
\lim_{n\rightarrow\infty} \frac{\ind{|\dG_n^*|>0}}{|\dT_{n-1}^*|}\bs{U}_{n-1}\bs{\wh{\sigma}}_n
&\!\!\!=\!\!\!&\lim_{n\rightarrow\infty} \frac{\ind{|\dG_n^*|>0}}{|\dT_{n-1}^*|}
\sum_{k \in \mathbb{T}_{n-1}}
\left( \begin{array}{c}
\wh{\eps}^2_{2k} +\wh{\eps}^2_{2k+1} \\
2X_k\wh{\eps}^2_{2k} \\
2X_k \wh{\eps}^2_{2k+1}\\
X^2_k( \wh{\eps}^2_{2k} + \wh{\eps}^2_{2k+1} )
\end{array}\right)=\bs{U}\bs{\sigma}\indnex,\\
 \lim_{n\rightarrow\infty} \frac{\ind{|\dG_n^*|>0}1}{|\dT_{n-1}^*|}\bs{V}_{n-1}\bs{\wh{\rho}}_n
 &\!\!\!=\!\!\!&\lim_{n\rightarrow\infty} \frac{\ind{|\dG_n^*|>0}1}{|\dT_{n-1}^*|} \sum_{k \in \mathbb{T}_{n-1}}
\left( \begin{array}{c}
\wh{\eps}_{2k} \wh{\eps}_{2k+1} \\
2X_k \wh{\eps}_{2k}  \wh{\eps}_{2k+1}\\
X^2_k\wh{\eps}_{2k}  \wh{\eps}_{2k+1} 
\end{array}\right)
=\bs{V}\bs{\rho}\indnex,
 \end{eqnarray*} 
 hence the result using Proposition~\ref{prop lim SUV}.
\hspace{\stretch{1}}$ \Box$

\section{Martingales and convergence rate}
\label{section martingales}
The aim of this section is to obtain sharper convergence results for our estimators, namely rates of convergence. The $P$-BMC approach does not allow this, therefore we now use martingale theory instead, as in \cite{BSG09,SGM11}. However, we cannot directly apply the results therein mainly because our noise sequence $(\eps_k=\veps_k+\eta_kX_{[k/2]})$ now contains the BAR process $(X_k)$ and thus does not satisfy the assumptions of \cite{BSG09,SGM11}.
\subsection{Martingales on binary trees}
\label{section meta lemme}
We start with a general result of convergence for martingales on a Galton-Watson binary tree, that we will make repeatedly use of in the following sections. Special cases of this result have already been proved and used in \cite{BSG09} and \cite{SGM11}. Note that in this binary tree context, we cannot use the standard asymptotic theory for vector martingales (see e.g. \cite{Duflo97}) because the number of data is roughly multiplied by $m$ at each generation.

\begin{Theorem}\label{meta martingale}
Let $(\bs{M}_n)$ be a $p$-dimensional $\dF^{\cO}$-martingale on the GW-binary tree $\dT^*$:
$\bs{M}_n=\sum_{\ell=1}^n\sum_{k\in\dG_{\ell}^*}\bs{W}_k$,
with $\bs{W}_k=(w_k^1,w_k^2,\ldots,w_k^p)^t$.  We make the following assumptions
\begin{description}
\item[(A.1)] $(\bs{M}_n)$ is square-integrable.
\end{description}
Let $<\bs{M}>_n=\sum_{\ell=0}^{n-1} \bs{\Gamma}_{\ell}$ be the predictable quadratic variation of $(\bs{M}_n)$, with
\begin{equation*}
\bs{ \Gamma}_{n}=\dE[\Delta\bs{M}_{n+1}\Delta\bs{M}_{n+1}^t\ |\ \cF^{\cO}_{n}].
\end{equation*}
\begin{description}
\item[(A.2)] $|\dT_{n-1}^*|^{-1}<\bs{M}>_n$ converges almost surely to a positive semi-definite  matrix $\bs{\Gamma}$ on $\overline{\cE}$.
\item[(A.3)] The $p\times p$ $\dF^{\cO}$-matrix martingale $(\bs{K}_n)$ defined by
\begin{equation*}
\bs{K}_n=\sum_{\ell=1}^n|\dT_{\ell}^*|^{-1}(\Delta\bs{M}_{\ell+1}\Delta\bs{M}_{\ell+1}^t-\dE[\Delta\bs{M}_{\ell+1}\Delta\bs{M}_{\ell+1}^t\ |\ \cF^{\cO}_{\ell}])
\end{equation*}
is square-integrable and its component-wise predictable quadratic variationes are $\cO(n)$ a.s.  on $\overline{\cE}$.
\end{description}
Let $(\bs{\Xi}_n)$ be a sequence of $p\times p$ invertible symmetric matrices such that
\begin{description}
\item[(A.4)] $|\dT_{n}^*|^{-1}\bs{\Xi}_n$ converges a.s. to an invertible matrix $\bs{\Xi}$  on $\overline{\cE}$;
\item[(A.5)] there exists a positive constant $\alpha$ such that on the non-extinction set $\overline{\cE}$ and for all $n$ and the following assumptions holds $\alpha(\bs{\Xi}_{n-1}^{-1}-\bs{\Xi}_n^{-1})\geq\bs{\Xi}_n^{-1}\bs{\Gamma}_n\bs{\Xi}_n^{-1}$ , in the sense that $\alpha(\bs{\Xi}_{n-1}^{-1}-\bs{\Xi}_n^{-1})-\bs{\Xi}_n^{-1}\bs{\Gamma}_n\bs{\Xi}_n^{-1}$ is a positive semi-definite  matrix.
\end{description}
Then $\bs{M}_n^t\bs{\Xi}_{n-1}^{-1}\bs{M}_n=\cO(n)$ and if $\bs{\Xi}$ is positive definite  $\|\bs{M}_n\|^2=\cO(nm^n)$ a.s.\\
If in addition, the entries of $(\bs{M}_n)$ satisfy
\begin{description}
\item[(A.6)] $\sup_n\dE[(m^{-n/2}\sum_{k\in\dG_{n}^*}w_k^i)^4\ |\ \cF^{\cO}_{n-1}]<\infty$ almost surely,
\end{description}
then for all $\delta>1/2$, $\|\bs{M}_n\|^2=o(n^{\delta}m^n)$ a.s. and
\begin{equation*}
\lim_{n\rightarrow\infty}\ind{|\dG_n^*|>0}\frac{1}{n}\sum_{\ell=1}^n\bs{M}_{\ell}^t\bs{\Xi}_{\ell-1}^{-1}\bs{M}_{\ell}=tr(\bs{\Gamma}\bs{\Xi}^{-1})\indnex\quad\textrm{a.s.}
\end{equation*}
\end{Theorem}

\noindent\textbf{Proof of the first part of Theorem~\ref{meta martingale}:} The result is obvious on the extinction set $\cE$. In the sequel, let us suppose that we are on the non-extinction set $\overline{\cE}$. For all $n\geq 1$, denote $\cV_{n}=\bs{M}_n^t\bs{\Xi}_{n-1}^{-1}\bs{M}_n$. The following equalities hold
\begin{eqnarray*}
\cV_{n+1}
&\!\!\!=\!\!\!& \bs{M}_{n+1}^t\bs{\Xi}_n^{-1}\bs{M}_{n+1}=(\bs{M}_n+\Delta \bs{M}_{n+1})^t\bs{\Xi}_n^{-1}(\bs{M}_n+\Delta \bs{M}_{n+1}),\\
&\!\!\!=\!\!\!&\cV_n-\bs{M}_n^{t}(\bs{\Xi}_{n-1}^{-1}-\bs{\Xi}_n^{-1})\bs{M}_n+
2\bs{M}_n^t\bs{\Xi}_n^{-1}\Delta \bs{M}_{n+1}\!+\!\Delta \bs{M}_{n+1}^t\bs{\Xi}_n^{-1}\Delta \bs{M}_{n+1},
\end{eqnarray*}
since $\bs{M}_n^t\bs{\Xi}_n^{-1}\Delta \bs{M}_{n+1}$ is scalar, and hence equal to its own transpose. By summing over the identity above, we obtain
\begin{equation}
\label{maindecomart}
\cV_{n+1}+\cA_n=\cV_1+\cB_{n+1}+\cW_{n+1},
\end{equation}
where
\begin{equation*}
\cA_n=\sum_{\ell=1}^n \bs{M}_\ell^{t}(\bs{\Xi}_{\ell-1}^{-1}-\bs{\Xi}_\ell^{-1})\bs{M}_\ell,
\quad
\cB_{n+1}=2\sum_{\ell=1}^n \bs{M}_\ell^t\bs{\Xi}_\ell^{-1}\Delta \bs{M}_{\ell+1},\quad
\cW_{n+1}=\sum_{\ell=1}^n \Delta \bs{M}_{\ell+1}^t\bs{\Xi}_\ell^{-1}\Delta \bs{M}_{\ell+1}.
\end{equation*}
The asymptotic behavior of the sequences $(\cW_n)$ and $(\cB_n)$ is given in the following lemmas.

\begin{Lemma}
\label{lemlimW}
Under assumptions \emph{(\textbf{A.1})} to \emph{(\textbf{A.4})}, the following almost sure convergence holds
\begin{equation*}
\lim_{n \rightarrow \infty}\ind{|\dG_n^*|>0}\frac{1}{n} \cW_n= \frac{m-1}{m}tr(\bs{\Gamma}\bs{\Xi}^{-1})\indnex
\hspace{1cm}\text{a.s.}
\end{equation*}
\end{Lemma}

\begin{Lemma}
\label{lemlimB}
Under assumptions \emph{(\textbf{A.1})} to \emph{(\textbf{A.5})}, the following asymptotic property holds
\begin{equation*}
\cB_{n+1}=o(n)\hspace{1cm}\textrm{a.s.}
\end{equation*}
\end{Lemma}
One then obtains
\begin{equation}
\label{cvgVA}
\lim_{n \rightarrow \infty}\ind{|\dG_n^*|>0}\frac{\cV_{n+1}+\cA_n}{n} = \frac{m-1}{m}tr(\bs{\Gamma}\bs{\Xi}^{-1})\indnex
\hspace{1cm}\text{a.s.}
\end{equation}
As $(\cA_n)$ is a sequence of positive real numbers, it follows that $\cV_{n+1}=\cO(n)$ a.s., which means that
$
\bs{M}_n^t\bs{\Xi}_{n-1}^{-1}\bs{M}_n=\cO(n)
$ a.s. 
As $\bs{\Xi}$ is positive definite, one obtains for large enough $n$, on the non-extinction set $\overline{\cE}$
\begin{equation*}
\|\bs{M}_n\|^2=\bs{M}_n^t\bs{M}_n\leq \frac{\bs{M}_n^t\bs{\Xi}_{n-1}^{-1}\bs{M}_n}{\lambda_{\textrm{min}}(\bs{\Xi}_{n-1}^{-1})},
\end{equation*}
where $\lambda_{\textrm{min}}(\bs{\Xi}_{n-1}^{-1})$ denotes the smallest eigenvalue of $\bs{\Xi}_{n-1}^{-1}$.
Finally, Assumption \textbf{(A.4)} yields
$
\|\bs{M}_n\|^2=\cO(nm^n)
$ a.s.,
which completes the proof of the first part of Theorem~\ref{meta martingale}.
\hspace{\stretch{1}}$\Box$\\

\noindent It remains to prove Lemmas~\ref{lemlimW} and \ref{lemlimB}.\\

\noindent\textbf{Proof of Lemma~\ref{lemlimW}:} First of all, we decompose $\cW_{n+1}=\cT_{n+1}+\cR_{n+1}$ with
\begin{equation*}
\cT_{n+1}=\sum_{\ell=1}^n \frac{\Delta \bs{M}_{\ell+1}^t\bs{\Xi}^{-1}\Delta \bs{M}_{\ell+1}}{|\dT_{\ell}^*|},\quad
\cR_{n+1}=\sum_{\ell=1}^n \frac{\Delta \bs{M}_{\ell+1}^t(|\dT_{\ell}^*|\bs{\Xi}_{\ell}^{-1}-\bs{\Xi}^{-1})\Delta \bs{M}_{\ell+1}}{|\dT_{\ell}^*|}.
\end{equation*}
We first prove that
$
\lim_{n \rightarrow \infty}\ind{|\dG_n^*|>0}\frac{1}{n} \cT_n= \frac{m-1}{m}tr(\bs{\Gamma}\bs{\Xi}^{-1})\indnex
$ a.s. 
As  $\cT_{n}$ is a scalar and the trace is commutative, we can rewrite $\cT_{n+1}=tr(\bs{H}_{n+1}\bs{\Xi}^{-1})$ where
\begin{equation*}
\bs{H}_{n+1}=\sum_{\ell=1}^{n} \frac{\Delta \bs{M}_{\ell+1}\Delta \bs{M}_{\ell+1}^t}{|\dT_{\ell}^*|}=\sum_{\ell=1}^{n} \frac{\bs{\Gamma}_{\ell}}{|\dT_{\ell}^*|}+\bs{K}_{n}.
\end{equation*}
On the one hand, by assumption \textbf{(A.3)}, one has $\bs{K}_n=o(n)$ a.s. on $\overline{\cE}$. On the other hand, Assumption \textbf{(A.2)} yields
\begin{equation*}
\ind{|\dG_{\ell}^*|>0} \frac{\bs{\Gamma}_{\ell}}{|\dT_{\ell}^*|}
=\ind{|\dG_{\ell}^*|>0}\left(\frac{<\bs{M}>_{\ell+1}}{|\dT_{\ell}^*|}-\frac{|\dT_{\ell-1}^*|}{|\dT_{\ell}^*|}\frac{<\bs{M}>_{\ell}}{|\dT_{\ell-1}^*|}
\right)
\end{equation*}
so that $H_n$ converges to $\left(\bs{\Gamma}-\bs{\Gamma}/m\right)\indnex=\bs{\Gamma}(m-1)/m$ a.s. as $\ell$ tends to infinity. Hence, Cesaro convergence yields 
\begin{equation}\label{cv Hn}
\lim_{n\rightarrow\infty}\ind{|\dG_n^*|>0}\frac{1}{n}\bs{H}_{n}=\frac{m-1}{m}\bs{\Gamma}\indnex \qquad\text{a.s.}
\end{equation} 
As a consequence, $\lim\ind{|\dG_n^*|>0}\cT_{n}/n=tr(\bs{\Gamma}\bs{\Xi}^{-1})(m-1)/m\indnex$ a.s.
We now turn to the asymptotic behavior of $\cR_{n+1}$. We know from Assumption \textbf{(A.4)} that $|\dT_{\ell}^*|\bs{\Xi}_{\ell}^{-1}-\bs{\Xi}^{-1}$ goes to $0$ as $\ell$ goes to infinity on $\overline{\cE}$. Thus, for all positive $\epsilon$, there exists $\ell_{\epsilon}$ such that if $\ell \ge \ell_{\epsilon}$,
\begin{equation*}
\ind{|\dG_\ell^*|>0}|\Delta \bs{M}_{\ell+1}^t(|\dT_{\ell}|\bs{\Xi}_{\ell}^{-1}-\bs{\Xi}^{-1})\Delta \bs{M}_{\ell+1}|
\leq 4\epsilon \Delta \bs{M}_{\ell+1}^t\Delta \bs{M}_{\ell+1}\ind{|\dG_\ell^*|>0}.
\end{equation*}
Hence, there exists some positive real number $c_{\epsilon}$ such that, for $n \ge \ell_{\epsilon}$,
\begin{eqnarray*}
\ind{|\dG_n^*|>0}|\cR_n| &\le& \ind{|\dG_n^*|>0} \Big(4\epsilon \sum_{\ell=\ell_{\epsilon}}^{n-1}\ind{|\dG_{\ell}^*|>0}\frac{\Delta \bs{M}_{\ell+1}^t \Delta \bs{M}_{\ell+1}}{|\dT^*_{\ell}|} + c_{\epsilon}\Big) \\
&\le&  \ind{|\dG_n^*|>0} \Big(4\epsilon\  tr(\bs{H}_n) + c_{\epsilon}\Big).
\end{eqnarray*}
This last inequality holding for any positive  $\epsilon$ and large enough $n$, the limit given by Equation~(\ref{cv Hn}) entails that
$
\lim_{n \rightarrow \infty}\ind{|\dG_n^*|>0}\frac{1}{n} \cR_n= 0
$ a.s.,
which completes the proof.
\hspace{\stretch{1}}$\Box$\\

\noindent\textbf{Proof of Lemma~\ref{lemlimB}:} Again, the result is obvious on the extinction set $\cE$. Suppose now we are on the non-extinction set $\overline{\cE}$. Recall that
\begin{equation*}
\cB_{n+1}=2\sum_{k=1}^n \bs{M}_k^t\bs{\Xi}_k^{-1}\Delta \bs{M}_{k+1}.
\end{equation*}
The process $(\cB_n)$ is a real-valued $\dF^{\cO}$-martingale. In addition, the following equality clearly holds
\begin{equation*}
\dE[\Delta\cB_{n+1}^2|\cF_n^{\cO}]=4\bs{M}_n^t\bs{\Xi}_n^{-1}\bs{\Gamma}_n\bs{\Xi}_n^{-1}\bs{M}_n
\hspace{1cm}\text{a.s.}
\end{equation*}
Assumption \textbf{(A.5)} then yields
\begin{equation*}
<\! \cB \!>_{n+1}\leq 4\alpha\sum_{k=1}^n\bs{M}_k^t(\bs{\Xi}_{k-1}^{-1}-\bs{\Xi}_k^{-1})\bs{M}_k=4\alpha\cA_n.
\hspace{1cm}\text{a.s.}
\end{equation*}
Hence, the law of large number for real martingales yields $\cB_n=o(\cA_n)$ a.s. Finally, we deduce from decomposition (\ref{maindecomart}) and Lemma~\ref{lemlimW} that
\begin{equation*}
\cV_{n+1}+ \cA_n=o(\cA_n)+\cO(n)
\hspace{1cm}\text{a.s.}
\end{equation*}
leading to $\cA_{n}=\cO(n)$ and $\cV_{n+1}=\cO(n)$ a.s. as both sequences are non-negative. This implies in turn that $\cB_n=o(n)$ a.s. completing the proof.
\hspace{\stretch{1}}$\Box$\\

\noindent\textbf{Proof of the second part of Theorem~\ref{meta martingale}:} Let us rewrite the entries $M_n^q$ of the martingale $\bs{M}_n$ as
$$M_n^q = \sum_{l=1}^n \underbrace{m^{\ell/2}}_{u_\ell} \underbrace{\frac{1}{m^{\ell/2}}\sum_{k \in \dG_l} w_k^q}_{\tilde{\varepsilon}_\ell^q}$$
Then one just has to apply Wei's lemma given in \cite[p 1672]{Wei87} to the martingale difference sequence $\tilde{\varepsilon}_\ell$ and $u_{\ell}=m^{\ell/2}$, and for the function $f(x)=(\log x)^{\delta}$ for $\delta>1/2$. Under assumption \textbf{(A.6)}, one obtains $M_n^q=o(m^{n/2}n^{\delta/2})$. As $P_n^q$ is the $q$-th entry of $\bs{M}_n$, one obtains $\|\bs{M}_n\|^2=o(n^{\delta}m^n)$ a.s. Now recall that $\cV_{n}=\bs{M}_n^t\bs{\Xi}_{n-1}^{-1}\bs{M}_n$, therefore, the following equality holds
\begin{equation*}
\ind{|\dG_n^*|>0}\cV_n= \ind{|\dG_n^*|>0}\bs{M}_n^t\bs{\Xi}_{n-1}^{-1}\bs{M}_n=o(n^{2\delta})\quad\text{a.s.},
\end{equation*}
for all $\delta>1/4$. In particular, for $\delta=1/2$, and we have the following order  $\cV_{n}=o(n)$. Lemmas~\ref{lemlimW} and \ref{lemlimB} then yield 
\begin{equation}\label{lim An}
\lim_{n\rightarrow\infty}\ind{|\dG_n^*|>0}\frac{\cA_n}{n} =\frac{m-1}{m}tr(\bs{\Gamma}\bs{\Xi}^{-1})\indnex\quad\text{a.s.}
\end{equation}
First of all, $\cA_n$ may be rewritten as
\begin{equation*}
\cA_n
=\sum_{\ell=1}^n \bs{M}_{\ell}^{t}(\bs{\Xi}_{\ell-1}^{-1}-\bs{\Xi}_{\ell}^{-1})\bs{M}_{\ell}
=\sum_{\ell=1}^n\bs{M}_{\ell}^{t}\bs{\Delta}_{\ell} \bs{\Xi}_{\ell-1}^{-1}\bs{M}_{\ell},
\end{equation*}
where $\bs{\Delta}_n=\rI_{4} - \bs{\Xi}_n^{-1} \bs{\Xi}_{n-1}$.
Thanks to Assumption \textbf{(A.4)} we know that
\begin{equation*}
\lim_{n\rightarrow\infty}\ind{|\dG_n^*|>0}\bs{\Delta}_n=\frac{m-1}{m}\rI_{4}\indnex\hspace{1cm}\text{a.s.}
\end{equation*}
Besides, Eq.~(\ref{lim An}) yields that $\ind{|\dG_n^*|>0}\cA_n \sim n (m-1)m^{-1}tr(\bs{\Gamma}\bs{\Xi}^{-1})\indnex$ a.s. Plugging these two results into the equality
\begin{equation*} 
\cA_n=\frac{m-1}{m}\sum_{\ell=1}^n\bs{M}_{\ell}^{t}\bs{\Xi}_{\ell-1}^{-1} \bs{M}_{\ell} 
+ \sum_{\ell=1}^n\bs{M}_{\ell}^{t}(\bs{\Delta}_{\ell}-\frac{m-1}{m}\rI_4)\bs{\Xi}_{\ell-1}^{-1}\bs{M}_{\ell}
\end{equation*}
gives that $\ind{|\dG_n^*|>0}\sum_{\ell=1}^n\bs{M}_{\ell}^{t}\bs{\Xi}_{\ell-1}^{-1} \bs{M}_{\ell} \sim \ind{|\dG_n^*|>0}\cA_n m(m-1)^{-1}$ a.s.
Thus one obtains
\begin{equation*}
\lim_{n\rightarrow\infty}\ind{|\dG_n^*|>0}\frac{1}{n}\sum_{\ell=1}^n\bs{M}_{\ell}^t\bs{\Xi}_{\ell-1}^{-1}\bs{M}_{\ell}=tr(\bs{\Gamma}\bs{\Xi}^{-1})\indnex\quad\textrm{a.s.}
\end{equation*}
which is the expected result.
\hspace{\stretch{1}}$\Box$

\subsection{Rate of convergence for $\bs{\wh{\theta}}_n$}
\label{section theta}
We apply Theorem~\ref{meta martingale} to a suitably chosen martingale. Recall that
\begin{equation} \label{diff}
\wh{\bs{\theta}}_n-\bs{\theta} 
= \bs{S}^{-1}_{n-1}
\sum_{k \in \dT_{n-1}}
\left(
\eps_{2k},
X_k\eps_{2k} ,
\eps_{2k+1},
X_k\eps_{2k+1}
\right)^t=\bs{S}^{-1}_{n-1}\bs{M}_n,
\end{equation}
where
$$\bs{M}_n= \sum_{k \in \dT_{n-1}} \left( 
\eps_{2k},
X_k\eps_{2k},
\eps_{2k+1},
X_k\eps_{2k+1}\right)^t.$$
Under assumptions \textbf{(H.1-3)}, for all
$n\geq 0$, $k \in \dG_n$, $\dE[\eps_{2k+i}|\cF_n^{\cO}]=\dE[X_k\eps_{2k+i}|\cF_n^{\cO}]=0$ and $(\bs{M}_n)$ is  a square-integrable $(\cF_n^{\cO})$-martingale, so that Assumption~\textbf{(A.1)} of Theorem~\ref{meta martingale} holds. Let us compute the predictable quadratic variation of $(\bs{M}_n)$
\begin{equation}
\label{dcroch}
\dE[\Delta \bs{M}_{n+1}\Delta \bs{M}_{n+1}^t|\cF_n^{\cO}]=\bs{\Gamma}_{n}=\sum_{k \in \dG_{n}}\bs{\gamma}_k
\otimes\left(
\begin{array}{cc}
1&X_k\\
X_k&X_k^2
\end{array}
\right),
\end{equation}
where
\begin{equation}
\label{gamk}
\bs{\gamma}_k=
\left(
\begin{array}{rr}
\delta_{2k}(\sigma^2_{\veps}+2X_k\rho_{00}+X_k^2\sigma^2_{\eta})&\delta_{2k}\delta_{2k+1}(\rho_{\veps}+2X_k\rho+X_k^2\rho_{\eta})\\
\delta_{2k}\delta_{2k+1}(\rho_{\veps}+2X_k\rho+X_k^2\rho_{\eta})&\delta_{2k+1}(\sigma^2_{\veps}+2X_k\rho_{11}+X_k^2\sigma^2_{\eta})
\end{array}
\right).
\end{equation} 
Thus the predictable quadratic variation of $(\bs{M}_n)$  is given by
\begin{equation*}
\cM_n= \sum_{\ell=0}^{n-1} \bs{\Gamma}_{\ell}=\sum_{k \in \dT_{n-1}}\bs{\gamma}_k
\otimes\left(
\begin{array}{cc}
1&X_k\\
X_k&X_k^2
\end{array}
\right).
\end{equation*}

\begin{Lemma}\label{lem lim Gamma}
Under assumptions \emph{\textbf{(H.1-5)}} and if $\kappa\geq 4$, the following convergence holds
\begin{equation*}
\lim_{n\rightarrow\infty}\ind{|\dG_n^*|>0}\frac{\cM_n}{|\dT_{n-1}^*|}=\bs{\Gamma}\indnex=
\left(\begin{array}{cc}
\bs{\Gamma}^0&\bs{\Gamma}^{01}\\
\bs{\Gamma}^{01}&\bs{\Gamma}^1
\end{array}\right)\indnex\quad\text{a.s.},
\end{equation*}
where $\bs{\Gamma}^0$, $\bs{\Gamma}^{01}$ and $\bs{\Gamma}^1$ are the $2\times 2$ matrices defined by
\begin{equation*}
\bs{\Gamma}^i=\left(\begin{array}{cc}
\sigma^2_{\veps}\ell_i(0)+2\rho_{ii}\ell_i(1)+\sigma^2_{\eta}\ell_i(2)&\sigma^2_{\veps}\ell_i(1)+2\rho_{ii}\ell_i(2)+\sigma^2_{\eta}\ell_i(3)\\
\sigma^2_{\veps}\ell_i(1)+2\rho_{ii}\ell_i(2)+\sigma^2_{\eta}\ell_i(3)&\sigma^2_{\veps}\ell_i(2)+2\rho_{ii}\ell_i(3)+\sigma^2_{\eta}\ell_i(4)
\end{array}\right),
\end{equation*}
and
\begin{equation*}
\bs{\Gamma}^{01}=\left(\begin{array}{cc}
\rho_{\veps}\ell_{01}(0)+2\rho\ell_{01}(1)+\rho_{\eta}\ell_{01}(2)&\rho_{\veps}\ell_{01}(1)+2\rho\ell_{01}(2)+\rho_{\eta}\ell_{01}(3)\\
\rho_{\veps}\ell_{01}(1)+2\rho\ell_{01}(2)+\rho_{\eta}\ell_{01}(3)&\rho_{\veps}\ell_{01}(2)+2\rho\ell_{01}(3)+\rho_{\eta}\ell_{01}(4)
\end{array}\right).
\end{equation*}
In addition, $\bs{\Gamma}$ is positive definite.
\end{Lemma}

\noindent\textbf{Proof} This is a direct consequence of Proposition~\ref{prop LGN}.
\hspace{\stretch{1}}$\Box$\\

\noindent Hence, Assumption \textbf{(A.2)} holds if $\kappa\geq 4$. The process $(\bs{K}_n)$ is clearly a square-integrable martingale if $\gamma\geq 2$. It is not difficult to check that its component-wise predictable quadratic variation is at worst of the order of
\begin{equation*}
\sum_{\ell=1}^n\frac{1}{|\dT_{\ell}^*|^2}\sum_{k\in\dG_{\ell}}\delta_{2k+i}X_k^8.
\end{equation*}
Proposition~\ref{prop LGN} ensures that $|\dT_n^*|^{-1}\sum_{k\in\dG_n}\delta_{2k+i}X_k^4$ converges almost surely on $\overline{\cE}$ provided $\kappa\geq 4$, it is therefore bounded by some constant $C$. As a result, its square is also bounded by $C^2$ and $|\dT_{\ell}^*|^{-2}\sum_{k\in\dG_{\ell}}\delta_{2k+i}X_k^8\leq C^2$ a.s. on $\overline{\cE}$. Finally, one obtains that 
\begin{equation*}
\ind{|\dG_n^*|>0}\sum_{\ell=1}^n\frac{1}{|\dT_{\ell}^*|^2}\sum_{k\in\dG_{\ell}}\delta_{2k+i}X_k^8\leq C^2n\ind{|\dG_n^*|>0},
\end{equation*}
so that Assumption \textbf{(A.3)} holds if $\kappa\geq 4$.

We now introduce a new sequence of matrices $\bs{\Sigma}_n$. They are defined as a standardized version of the predictable quadratic variation of $(\bs{M}_n)$, with the variance coefficients $\sigma^2_{\veps}$ and $\sigma^2_{\eta}$ set to $1$ and all the covariance coefficients $\rho_{\veps}$,  $\rho_{\eta}$, $\rho_{ij}$ set to $0$, namely
\begin{equation*}
\bs{\Sigma}_n=\sum_{\ell=1}^n\bs{\Phi}_{\ell}\bs{\Phi}_{\ell}^t=\sum_{k\in\dT_n}(1+X_k^2)
\left(\begin{array}{cc}
\delta_{2k}&0\\
0&\delta_{2k+1}
\end{array}\right)
\otimes
\left(\begin{array}{cc}
1&X_k\\
X_k&X_k^2
\end{array}\right),
\end{equation*}
where $\bs{\Phi}_{n}$ is the $4\times 2^n$ matrix of the collection of the $4\times 1$ vectors $(1+X^2_{k})^{1/2}(\delta_{2k}, \delta_{2k}X_k,\delta_{2k+1}, \delta_{2k+1}X_k)^t$ for $k\in\dG_n=\{2^n, 2^n+1,\ldots,2^{n+1}-1\}$
\begin{equation*}
\bs{\Phi}_{n}=\left(
\begin{array}{rrr}
\delta_{2(2^n)}\sqrt{1+X^2_{2^n}}
&\cdots&\delta_{2(2^{n+1}-1)}\sqrt{1+X^2_{2^{n+1}-1}}\\
\delta_{2(2^n)}X_{2^n}\sqrt{1+X^2_{2^n}}
&\cdots&\delta_{2(2^{n+1}-1)}X_{2^{n+1}-1}\sqrt{1+X^2_{2^{n+1}-1}}\\
\delta_{2(2^n)+1}\sqrt{1+X^2_{2^n}}
&\cdots&\delta_{2(2^{n+1}-1)+1}\sqrt{1+X^2_{2^{n+1}-1}}\\
\delta_{2(2^n)+1}X_{2^n}\sqrt{1+X^2_{2^n}}
&\cdots&\delta_{2(2^{n+1}-1)+1}X_{2^{n+1}-1}\sqrt{1+X^2_{2^{n+1}-1}}
\end{array}
\right).
\end{equation*}
Note that $\bs{\Phi}_{\ell}\bs{\Phi}_{\ell}^t$ and hence $\bs{\Sigma}_n$ is positive definite as soon as the $X_k$ are not constant. The next result is again a direct consequence of Proposition~\ref{prop LGN}.

\begin{Lemma}\label{lem lim Sigma}
Under assumptions \emph{\textbf{(H.1-5)}} and if $\kappa\geq 4$, the following convergence holds
\begin{equation*}
\lim_{n\rightarrow\infty}\ind{|\dG_n^*|>0}\frac{\bs{\Sigma}_n}{|\dT_{n}^*|}=\bs{\Sigma}\indnex=
\left(\begin{array}{cc}
\bs{\Sigma}^0&0\\
0&\bs{\Sigma}^1
\end{array}\right)\indnex\quad\text{a.s.},
\end{equation*}
where $\bs{\Sigma}^i$ is the $2\times 2$ matrix
$
\bs{\Sigma}^i=\left(\begin{array}{cc}
\ell_i(0)+\ell_i(2)&\ell_i(1)+\ell_i(3)\\
\ell_i(1)+\ell_i(3)&\ell_i(2)+\ell_i(4)
\end{array}\right).
$
In addition, $\bs{\Sigma}$ is positive definite.
\end{Lemma}
As a result, Assumption \textbf{(A.4)} also holds if $\kappa\geq 4$, and we now turn to Assumption \textbf{(A.5)}.

\begin{Lemma}\label{lem A5}
Under assumptions \emph{\textbf{(H.1-5)}}, for all $\alpha>\max\{2\sigma^2_{\veps},2\sigma^2_{\eta}, \mu^0,\mu^1,\nu\}$ and for all $n$, the following inequality holds $\bs{\Sigma}_n^{-1}\bs{\Gamma}_n\bs{\Sigma}_n^{-1}\leq\alpha(\bs{\Sigma}_{n-1}^{-1}-\bs{\Sigma}_n^{-1})$, where
\begin{equation*}
\mu^i=\frac{1}{2}\left(\sigma^2_{\veps}+\sigma^2_{\eta}+\left({(\sigma^2_{\veps}-\sigma^2_{\eta})^2+4\rho_{ii}^2}\right)^{\frac{1}{2}}\right),
\ \nu=\sigma^2_{\veps}+\sigma^2_{\eta}+\left({(\sigma^2_{\veps}-\sigma^2_{\eta})^2+(\rho_{00}+\rho_{11})^2}\right)^{\frac{1}{2}}.
\end{equation*}
\end{Lemma}

\noindent\textbf{Proof} We first prove that for all such $\alpha$, $\bs{\Gamma}_n\leq \alpha\bs{\Phi}_{n}\bs{\Phi}_{n}^t$ holds. For all $k\in\dG_n$, let $D_k^0=\alpha(1+X_k^2)-(\sigma^2_{\veps}+2X_k\rho_{00}+X_k^2\sigma^2_{\eta})$, $D_k^1=\alpha(1+X_k^2)-(\sigma^2_{\veps}+2X_k\rho_{11}+X_k^2\sigma^2_{\eta})$ and $D_k^{01}=\rho^2_{\veps}+2X_k\rho+X_k^2\rho^2_{\eta}$ be the coefficients of $\alpha\bs{\Phi}_{n}\bs{\Phi}_{n}^t-\bs{\Gamma}_n$ up to the sum over $\dG_{n}$. We first need to prove that $D_k^i>0$ for all $k$. One can rewrite 
$
D_k^i=\alpha-\sigma^2_{\veps}-2\rho_{ii}X_k+(\alpha-\sigma^2_{\eta})X_k^2,
$
so that it is sufficient to prove that this second order polynomial in $X_k$ has no real root, as both its terms of degree $0$ and $2$ are positive by assumption. Its discriminant in a function of $\alpha$ given by
\begin{equation*}
\Delta(\alpha)=-(\alpha^2-\alpha(\sigma^2_{\veps}+\sigma^2_{\eta})+\sigma^2_{\veps}\sigma^2_{\eta}-\rho_{ii}^2).
\end{equation*}
This discriminant $\Delta(\alpha)$ is again a second order polynomial in $\alpha$. Therefore it is negative as soon as $\alpha$ is larger than its largest root $\mu^i$.
Second, we want to prove that $(D_k^{01})^2\leq D_k^0D_k^1$. The Cauchy-Schwarz inequality yields
\begin{equation*}
(D_k^{01})^2\leq (\sigma^2_{\veps}+2X_k\rho_{00}+X_k^2\sigma^2_{\eta})(\sigma^2_{\veps}+2X_k\rho_{11}+X_k^2\sigma^2_{\eta}).
\end{equation*}
Hence, one just has to check that
\begin{equation*}
(\sigma^2_{\veps}+2X_k\rho_{00}+X_k^2\sigma^2_{\eta})(\sigma^2_{\veps}+2X_k\rho_{11}+X_k^2\sigma^2_{\eta})\leq D_k^0D_k^1,
\end{equation*}
which boils down to proving that the second order polynomial
\begin{equation*}
\alpha-2\sigma^2_{\veps}-2(\rho_{00}+\rho_{11})X_k+(\alpha-2\sigma^2_{\eta})X_k^2
\end{equation*}
is non negative. Similar arguments as above yield that the preceding polynomial  is nonnegative as soon as $\alpha>\max\{2\sigma^2_{\veps},2\sigma^2_{\eta}, \nu\}$. Thus, for $u=(u_1,u_2,u_3,u_4)^t\in\dR^4$ the following lower bound holds
\begin{eqnarray*}
{u^t(\alpha\bs{\Phi}_{n}\bs{\Phi}_{n}^t- \bs{\Gamma}_n)u}
&=&\sum_{k\in\dG_n}\left(u_1\delta_{2k}(D_k^0)^{1/2}+u_2\delta_{2k}X_k(D_k^0)^{1/2}-u_3\delta_{2k+1}\frac{D_k}{(D_k^0)^{1/2}}-u_4\delta_{2k+1}X_k\frac{D_k}{(D_k^0)^{1/2}}\right)^2\\
&&+\left(
u_3\delta_{2k+1}\Big(D_k^1-\frac{(D_k^{01})^2}{D_k^0}\Big)^{1/2}+u_4\delta_{2k+1}X_k\Big(D_k^1-\frac{(D_k^{01})^2}{D_k^0}\Big)^{1/2}\right)^2\\
&\geq&0,
\end{eqnarray*}
hence $\bs{\Gamma}_n\leq \alpha\bs{\Phi}_{n}\bs{\Phi}_{n}^t$. To obtain the expected result, we use Riccati equation (see e.g. Lemma B.1 of \cite{BSG09}) and the definition of $\bs{\Sigma}_n$ to obtain
\begin{equation}\label{eq Riccati}
\bs{\Sigma}_n^{-1}=\bs{\Sigma}_{n-1}^{-1}-\bs{\Sigma}_{n-1}^{-1}\bs{\Phi}_{n}(\bs{\rI}_{2^n}+\bs{l}_n)^{-1}\bs{\Phi}_n^t\bs{\Sigma}_{n-1}^{-1},
\end{equation}
where $\bs{l}_n=\bs{\Phi}_n^t\bs{\Sigma}_{n-1}^{-1}\bs{\Phi}_n$.
By multiplying both sides by $\bs{\Phi}_n$, we obtain
\begin{eqnarray*}
\bs{\Sigma}_n^{-1}\bs{\Phi}_n &=& \bs{\Sigma}_{n-1}^{-1}\bs{\Phi}_n-\bs{\Sigma}_{n-1}^{-1}\bs{\Phi}_{n}(\bs{\rI}_{2^n}+\bs{l}_n)^{-1}\bs{l}_n,\\
&=& \bs{\Sigma}_{n-1}^{-1}\bs{\Phi}_n-\bs{\Sigma}_{n-1}^{-1}\bs{\Phi}_{n}(\bs{\rI}_{2^n}+\bs{l}_n)^{-1}(\bs{\rI}_{2^n}+\bs{l}_n-\bs{\rI}_{2^n}),\\
&=&\bs{\Sigma}_{n-1}^{-1}\bs{\Phi}_{n}(\bs{\rI}_{2^n}+\bs{l}_n)^{-1}.
\end{eqnarray*}
In particular, as $\bs{l}_n$ is positive definite, one obtains
\begin{equation*}
\bs{\Sigma}_n^{-1}\bs{\Phi}_n(\bs{\rI}_{2^n}+\bs{l}_n)^{1/2}=\bs{\Sigma}_{n-1}^{-1}\bs{\Phi}_{n}(\bs{\rI}_{2^n}+\bs{l}_n)^{-1/2}.
\end{equation*}
Taking the square of each side of the above equation then yields
\begin{equation*}
\bs{\Sigma}_n^{-1}\bs{\Phi}_n(\bs{\rI}_{2^n}+\bs{l}_n)\bs{\Phi}_n^t\bs{\Sigma}_n^{-1}=\bs{\Sigma}_{n-1}^{-1}\bs{\Phi}_{n}(\bs{\rI}_{2^n}+\bs{l}_n)^{-1}\bs{\Phi}_n^t\bs{\Sigma}_{n-1}^{-1}=\bs{\Sigma}_{n-1}^{-1}-\bs{\Sigma}_n^{-1},
\end{equation*}
by Equation~(\ref{eq Riccati}). As $\bs{\rI}_{2^n}+\bs{l}_n\geq \bs{\rI}_{2^n}$ in the sense of positive semi-definite  matrices, one obtains
$
\bs{\Sigma}_{n-1}^{-1}-\bs{\Sigma}_n^{-1}\geq\bs{\Sigma}_n^{-1}\bs{\Phi}_n\bs{\Phi}_n^t\bs{\Sigma}_n^{-1}.
$
This inequality together with $\bs{\Gamma}_n\leq \alpha\bs{\Phi}_{n}\bs{\Phi}_{n}^t$ yield the result.
\hspace{\stretch{1}}$\Box$

\begin{Lemma}\label{lem A5}
Under assumptions \emph{\textbf{(H.1-5)}} and if $\kappa\geq 4$, for $i\in\{0,1\}$ and $q\in\{0,1\}$, one obtains
\begin{equation*}
\sup_n\Big\{m^{-2n}\dE\big[\big(\sum_{k\in\dG_n^*}X_k^q\epsilon_{2k+i}\big)^4\ |\ \cF_n^{\cO}\big]\Big\}<\infty\qquad a.s.
\end{equation*}
\end{Lemma}

\noindent\textbf{Proof} The following inequality is easily proved
\begin{eqnarray*}
{m^{-2n}\dE\big[\big(\sum_{k\in\dG_n^*}X_k^q\epsilon_{2k+i}\big)^4\ |\ \cF_n^{\cO}\big]}
&\leq& C\Big(\frac{1}{m^n}\sum_{k\in\dG_n^*}\delta_{2k+i}X_k^{2q}(1+X_k+X_k^2)\Big)^2\\
&&+C\frac{1}{m^{2n}}\sum_{k\in\dG_n^*}\delta_{2k+i}X_k^{4q}(1+X_k+X_k^2+X^3+X^4),
\end{eqnarray*}
where $C$ is a constant depending only on the moments of $(\veps_2,\eta_2,\veps_3,\eta_3)$ up to order 4. The result follows from Proposition~\ref{prop LGN}.
\hspace{\stretch{1}}$\Box$\\

\noindent We have now proved that Assumptions \textbf{(A.1-6)} of Theorem~\ref{meta martingale} hold for the martingale $(\bs{M}_n)$ and the sequence of positive definite matrices $(\bs{\Xi}_n)=(\bs{\Sigma}_n)$, thus we obtain the following result.

\begin{Proposition}\label{prop vitesse M}
Under assumptions \emph{\textbf{(H.1-5)}} and if $\kappa\geq 4$, one obtains 
\begin{equation*}
\bs{M}_n^t\bs{\Sigma}_{n-1}^{-1}\bs{M}_n=\cO(n),\quad\textrm{and}\quad \|\bs{M}_n\|^2=\cO(nm^n)\quad\textrm{a.s.}
\end{equation*}
In addition, for all $\delta>1/2$, $\|\bs{M}_n\|^2=o(n^{\delta}m^n)$ a.s. and 
\begin{equation*}
\lim_{n\rightarrow\infty}\ind{|\dG_n^*|>0}\frac{1}{n}\sum_{\ell=1}^n\bs{M}_{\ell}^t\bs{\Sigma}_{\ell-1}^{-1}\bs{M}_{\ell}=tr(\bs{\Gamma}\bs{\Sigma}^{-1})\indnex,\quad\textrm{a.s.}
\end{equation*}
\end{Proposition}

Now recall that $\wh{\bs{\theta}}_n-\bs{\theta} = \bs{S}^{-1}_{n-1}\bs{M}_n$. One then readily obtains Theorem~\ref{th cv theta chap}.\\

\noindent\textbf{Proof of Theorem~\ref{th cv theta chap}} 
As  $\wh{\bs{\theta}}_n-\bs{\theta} = \bs{S}^{-1}_{n-1}\bs{M}_n$, one obtains
\begin{equation*}
\|\wh{\bs{\theta}}_n-\bs{\theta}\|^2=\bs{M}_n^t\bs{S}^{-2}_{n-1}\bs{M}_n\leq{\|\bs{M}_n\|^2}{\lambda_{\textrm{max}}(\bs{S}^{-2}_{n-1})}.
\end{equation*}where $\lambda_{\textrm{max}}(\bs{S}^{-2}_{n-1})$ denotes the highest eigenvalue of matrix $\bs{S}^{-2}_{n-1}$.
We use Proposition~\ref{prop lim SUV} to conclude that $\|\wh{\bs{\theta}}_n-\bs{\theta} \|^2=o(n^{\delta}m^{-n})$ a.s. For the quadratic strong law, Proposition~\ref{prop vitesse M} yields
\begin{equation*}
\lim_{n\rightarrow\infty}\ind{|\dG_n^*|>0}\frac{1}{n}\sum_{\ell=1}^n(\wh{\bs{\theta}}_{\ell}-\bs{\theta})^t\bs{S}_{\ell-1}\bs{\Sigma}_{\ell-1}^{-1}\bs{S}_{\ell-1}(\wh{\bs{\theta}}_{\ell}-\bs{\theta})=tr(\bs{\Gamma}\bs{\Sigma}^{-1})\indnex\qquad\textrm{a.s.}
\end{equation*}
and the result is obtained by using Proposition~\ref{prop lim SUV} and Lemma~\ref{lem lim Sigma}. A similar argument as in the proof of Lemma~\ref{lemlimW} is used to replace $\bs{S}_{\ell-1}\bs{\Sigma}_{\ell-1}^{-1}\bs{S}_{\ell-1}$ by its equivalent $|\dT_{\ell-1}^*|^{-1}\bs{S}\bs{\Sigma}^{-1}\bs{S}$.
\hspace{\stretch{1}}$\Box$

\subsection{Rate of convergence for $\bs{\wh{\sigma}}_n$}
\label{section sigma}
We proceed in two steps. Recall that 
\begin{equation*}
\wh{{\bs{\sigma}}}_n\ =\ 
\bs{U}^{-1}_{n-1}
\sum_{k \in \mathbb{T}_{n-1}}
\left( 
\wh{\eps}^2_{2k} +\wh{\eps}^2_{2k+1},
2X_k\wh{\eps}^2_{2k},
2X_k \wh{\eps}^2_{2k+1},
X^2_k( \wh{\eps}^2_{2k} + \wh{\eps}^2_{2k+1} )
\right)^t,
\end{equation*}
is our estimator of $\bs{\sigma}=(\sigma^2_{\veps}, \rho_{00}, \rho_{11}, \sigma^2_{\eta})^t$, and
\begin{equation*}
{{\bs{\sigma}}}_n\ =\ 
\bs{U}^{-1}_{n-1}
\sum_{k \in \mathbb{T}_{n-1}}
\left( 
{\eps}^2_{2k} +{\eps}^2_{2k+1} ,
2X_k{\eps}^2_{2k} ,
2X_k{\eps}^2_{2k+1},
X^2_k( {\eps}^2_{2k} + {\eps}^2_{2k+1} )
\right)^t.
\end{equation*}
Our first step is to prove the convergence of ${{\bs{\sigma}}}_n$ to ${{\bs{\sigma}}}$. The second step is the convergence of $\wh{{\bs{\sigma}}}_n-{{\bs{\sigma}}}_n$ with a convergence rate.

\subsubsection{Convergence of ${{\bs{\sigma}}}_n$ }
\label{section cv sigma n}
The convergence of  ${{\bs{\sigma}}}_n$ to ${{\bs{\sigma}}}$ is directly obtained using the usual law of large numbers for square-integrable vector-valued martingales. Note that one could obtain a convergence rate under stronger moment assumptions using Theorem~\ref{meta martingale}.

\begin{Lemma}\label{lem sigman}
Under assumptions \emph{\textbf{(H.1-5)}} and if $\kappa\geq 8$, the following convergence holds
\begin{equation*}
\lim_{n\rightarrow\infty} \ind{|\dG_n^*|>0}{{\bs{\sigma}}}_n={{\bs{\sigma}}}\indnex\qquad a.s.
\end{equation*}
\end{Lemma}

\noindent\textbf{Proof :} Set 
\begin{equation*}
\bs{M}^{\bs{\sigma}}_n = \bs{U}_{n-1}({{\bs{\sigma}}}_n-{{\bs{\sigma}}})=\sum_{\ell=1}^{n-1}\sum_{k \in \mathbb{G}_{\ell}}
\left( \begin{array}{c}
{\eps}^2_{2k} +{\eps}^2_{2k+1} -\dE[{\eps}^2_{2k} +{\eps}^2_{2k+1}\ |\ \cF_{\ell}^{\cO}]\\
2X_k({\eps}^2_{2k} -\dE[{\eps}^2_{2k}\ |\ \cF_{\ell}^{\cO}])\\
2X_k({\eps}^2_{2k+1}-\dE[{\eps}^2_{2k+1}\ |\ \cF_{\ell}^{\cO}])\\
X^2_k( {\eps}^2_{2k} + {\eps}^2_{2k+1}-\dE[{\eps}^2_{2k} +{\eps}^2_{2k+1}\ |\ \cF_{\ell}^{\cO}] )
\end{array}\right).
\end{equation*}
Hence, $(\bs{M}^{\bs{\sigma}}_n)$ is a square-integrable $(\cF_n^{\cO})$-martingale. One can compute all the entries of its predictable quadratic variation and prove that they all equal a constant (depending on the moments of $(\veps_2,\eta_2,\veps_3,\eta_3)$ up to order 4) multiplied by $\sum\delta_{2k+i}X_k^q$ or $\sum\delta_{2k}\delta_{2k+1}X_k^q$ with $q\leq 8$. Hence, the laws of large numbers given in Proposition~\ref{prop LGN} ensure that  $m^{-n}<\bs{M}^{\bs{\sigma}}>_n$ converges almost surely to a constant matrix on the non-extinction set $\overline{\cE}$. The standard law of large numbers for square-integrable martingales then implies that $(\bs{M}^{\bs{\sigma}}_n)=o(m^n)$ a.s. 
Besides, $m^{-n}\bs{U}_n$ also converges to a fixed matrix on the non-extinction set $\overline{\cE}$ by Proposition~\ref{prop lim SUV}. Therefore $\ind{|\dG_n^*|>0}({{\bs{\sigma}}}_n-{{\bs{\sigma}}})= \ind{|\dG_n^*|>0}\bs{U}_{n-1}^{-1}\bs{M}^{\bs{\sigma}}_n$ tends to $0$ a.s. when $n$ tends to infinity.
\hspace{\stretch{1}}$\Box$

\subsubsection{Convergence of $\wh{{\bs{\sigma}}}_n-{{\bs{\sigma}}}_n$ }
\label{section cv sigma hat n}
We now turn to the convergence of  $\wh{{\bs{\sigma}}}_n-{{\bs{\sigma}}}_n$. One can rewrite $\bs{U}_{n-1}(\wh{{\bs{\sigma}}}_n-{{\bs{\sigma}}}_n)$ as
\begin{equation}\label{decomp sigma}
\bs{U}_{n-1}(\wh{{\bs{\sigma}}}_n-{{\bs{\sigma}}}_n)=\bs{P}_n^{\bs{\sigma}}+2\bs{R}_n^{\bs{\sigma}},
\end{equation}
with
\begin{equation*}
\bs{P}_n^{\bs{\sigma}}=\sum_{k \in \mathbb{T}_{n-1}}
\left( \begin{array}{c}
(\wh{\eps}_{2k}-{\eps}_{2k})^2 +(\wh{\eps}_{2k+1} -{\eps}_{2k+1} )^2\\
2X_k(\wh{\eps}_{2k}-{\eps}_{2k})^2 \\
2X_k(\wh{\eps}_{2k+1} -{\eps}_{2k+1} )^2\\
X^2_k\big( (\wh{\eps}_{2k}-{\eps}_{2k})^2 +(\wh{\eps}_{2k+1} -{\eps}_{2k+1} )^2 \big)
\end{array}\right),
\end{equation*}
and
\begin{equation*}
\bs{R}_n^{\bs{\sigma}}=\sum_{k \in \mathbb{T}_{n-1}}
\left( \begin{array}{c}
{\eps}_{2k}(\wh{\eps}_{2k}-{\eps}_{2k}) +{\eps}_{2k+1}(\wh{\eps}_{2k+1} -{\eps}_{2k+1} )\\
2X_k{\eps}_{2k}(\wh{\eps}_{2k}-{\eps}_{2k}) \\
2X_k{\eps}_{2k+1}(\wh{\eps}_{2k+1} -{\eps}_{2k+1} )\\
X^2_k\big( {\eps}_{2k}(\wh{\eps}_{2k}-{\eps}_{2k}) +{\eps}_{2k+1}(\wh{\eps}_{2k+1} -{\eps}_{2k+1} ) \big)
\end{array}\right).
\end{equation*}
We are going to study separately the asymptotic properties of $\bs{P}_n^{\bs{\sigma}}$ and $\bs{R}_n^{\bs{\sigma}}$.

\begin{Lemma}\label{lem eps2-0}
Under assumptions \emph{\textbf{(H.1-5)}} and if $\kappa\geq 4$, one obtains
\begin{equation*}
\lim_{n\rightarrow\infty}\ind{|\dG_n^*|>0}\frac{1}{n}\sum_{k \in \dT_{n}} (\wh{\eps}_{2k}-{\eps}_{2k})^2=q_0(0)\indnex=(m-1)tr(\bs{\Gamma}^0(\bs{S}^0)^{-1})\indnex\qquad a.s.
\end{equation*}
where $\bs{\Gamma}^0$ is defined in Lemma~\ref{lem lim Gamma} and $\bs{S}^0$ in Proposition~\ref{prop lim SUV}.
\end{Lemma}

\noindent\textbf{Proof :} We are going to apply Theorem~\ref{meta martingale} to the first two entries of the martingale $(\bs{M}_n)$. Indeed, let $\bs{M}_n^0$ be the $2$-component vector corresponding to the first two entries of $\bs{M}_n$
\begin{equation*}
\bs{M}_n^0=\sum_{k \in \dT_{n-1}} \left( \begin{array}{cccc}
\eps_{2k}  \\
X_k\eps_{2k}
\end{array}\right).
\end{equation*}
Let $\bs{\theta}^0=(a,b)^t$, $\wh{\bs{\theta}}^0_n=(\wh{a}_n,\wh{b}_n)^t$. Clearly, one has $(\wh{\bs{\theta}}^0_n-{\bs{\theta}}^0)=(\bs{S}_{n-1}^0)^{-1}\bs{M}_n^0$, therefore one obtains
\begin{eqnarray*}
\sum_{k \in \dG_{n}} (\wh{\eps}_{2k}-{\eps}_{2k})^2&=&\sum_{k \in \dG_{n}} \delta_{2k}(\wh{\bs{\theta}}^0_n-{\bs{\theta}}^0)^t
\left(\begin{array}{cc}
1&X_k\\
X_k&X_k^2
\end{array}\right)
(\wh{\bs{\theta}}^0_n-{\bs{\theta}}^0)\\
&=&(\bs{M}_n^0)^t(\bs{S}_{n-1}^0)^{-1}(\bs{S}_{n}^0-\bs{S}_{n-1}^0)(\bs{S}_{n-1}^0)^{-1}\bs{M}_n^0.
\end{eqnarray*}
Set also $\bs{\Sigma}_n^0$ to be the $2\times2$ matrix defined by
\begin{equation*}
\bs{\Sigma}_n^0=\sum_{k\in\dT_n}(1+X_k^2)
\delta_{2k}
\left(\begin{array}{cc}
1&X_k\\
X_k&X_k^2
\end{array}\right).
\end{equation*}
Thus, Proposition~\ref{prop lim SUV} and Lemma~\ref{lem lim Sigma} yield
\begin{equation*}\label{def delta0}
\lim_{n\rightarrow\infty}\ind{|\dG_n^*|>0}(\bs{\Sigma}_n^0)^{1/2}(\bs{S}_{n}^0)^{-1}(\bs{S}_{n+1}^0-\bs{S}_{n}^0)(\bs{S}_{n}^0)^{-1}(\bs{\Sigma}_n^0)^{1/2}=\bs{\Delta^0}\indnex\quad a.s.
\end{equation*}
where $\bs{\Delta^0}=(m-1)(\bs{\Sigma}^0)^{1/2}(\bs{S}^0)^{-1}(\bs{\Sigma}^0)^{1/2}$. Note that $\bs{\Delta^0}$ is positive definite and the matrices $\bs{\Sigma}^0_n$, $\bs{S}^0_n$, $\bs{\Sigma}^0$, $\bs{S}^0$ and $\bs{\Delta^0}$ commute. We now use Theorem~\ref{meta martingale} for the martingale $(\bs{M}_n^0)$, with the sequence  $\big(\Xi_n=(\bs{\Delta}^0)^{-1/2}\bs{\Sigma}^0_n(\bs{\Delta}^0)^{-1/2}\big)$. As $\bs{\Delta}^0$ is a fixed positive definite  matrix, it is clear that all the assumptions of Theorem~\ref{meta martingale} hold, as in Section~\ref{section theta}. Thus one obtains the a.s. limit
\begin{eqnarray*}
{\lim_{n\rightarrow\infty}\ind{|\dG_n^*|>0}\frac{1}{n}\sum_{\ell=1}^n(\bs{M}_{\ell}^0)^t((\bs{\Delta}^0)^{-1/2}\bs{\Sigma}_{\ell-1}^0(\bs{\Delta}^0)^{-1/2})^{-1}\bs{M}_{\ell}^0}
&=&tr(\bs{\Gamma}^0(\bs{\Sigma}^0)^{-1}\bs{\Delta}^0)\indnex\ =\ (m-1)tr(\bs{\Gamma}^0(\bs{S}^0)^{-1}).
\end{eqnarray*}
Finally, a similar argument as in the proof of Lemma~\ref{lemlimW} is used to replace $\bs{\Delta}^0$ by its asymptotic equivalent $(\bs{\Sigma}_n^0)^{1/2}(\bs{S}_{n}^0)^{-1}(\bs{S}_{n+1}^0-\bs{S}_{n}^0)(\bs{S}_{n}^0)^{-1}(\bs{\Sigma}_n^0)^{1/2}$ to obtain
\begin{eqnarray*}
{\lim_{n\rightarrow\infty}\ind{|\dG_n^*|>0}\frac{1}{n}\sum_{k \in \dT_{n}} (\wh{\eps}_{2k}-{\eps}_{2k})^2}
&=&\lim_{n\rightarrow\infty}\ind{|\dG_n^*|>0}\frac{1}{n}\sum_{\ell=1}^n (\bs{M}_{\ell}^0)^t(\bs{\Sigma}_{\ell-1}^0)^{-1/2}\bs{\Delta}^0(\bs{\Sigma}_{\ell-1}^0)^{-1/2}\bs{M}_{\ell}^0\nonumber\\
&=&\lim_{n\rightarrow\infty}\ind{|\dG_n^*|>0}\frac{1}{n}\sum_{\ell=1}^n (\bs{M}_{\ell}^0)^t\big((\bs{\Delta}^0)^{-1/2}\bs{\Sigma}_{\ell-1}^0(\bs{\Delta}^0)^{-1/2}\big)^{-1}\bs{M}_{\ell}^0,
\end{eqnarray*}
hence the result.
\hspace{\stretch{1}}$\Box$\\

\noindent A similar proof yields the following results for odd indices.

\begin{Lemma}\label{lem eps2-1}
Under assumptions \emph{\textbf{(H.1-5)}} and if $\kappa\geq 4$, the following convergence holds
\begin{equation*}
\lim_{n\rightarrow\infty}\ind{|\dG_n^*|>0}\frac{1}{n}\sum_{k \in \dT_{n}} (\wh{\eps}_{2k+1}-{\eps}_{2k+1})^2=q_1(0)\indnex=(m-1)tr(\bs{\Gamma}^1(\bs{S}^1)^{-1})\indnex\quad a.s.
\end{equation*}
where $\bs{\Gamma}^1$ is defined in Lemma~\ref{lem lim Gamma} and $\bs{S}^1$ in Proposition~\ref{prop lim SUV}.
\end{Lemma}
The proof of Lemma~\ref{lem eps2-0} can also be adapted to obtain the two following results.

\begin{Lemma}\label{lem esp2-2}
Under assumptions \emph{\textbf{(H.1-5)}} and if $\kappa\geq 4$, for all $i\in\{0,1\}$,  the almost sure convergence holds
\begin{equation*}
\lim_{n\rightarrow\infty}\ind{|\dG_n^*|>0}\frac{1}{n}\sum_{k \in \dT_{n}} X_k(\wh{\eps}_{2k+i}-{\eps}_{2k+i})^2=q_i(1)\indnex=(m-1)tr(\bs{\Gamma}^i(\bs{S}^i)^{-2}\bs{T}^i)\indnex,
\end{equation*}
where $\bs{T}^i$ is the $2\times 2$ matrix defined by
$
\bs{T}^i=\left(
\begin{array}{cc}
\ell_i(1)&\ell_i(2)\\
\ell_i(2)&\ell_i(3)
\end{array}\right).
$
\end{Lemma}

\noindent\textbf{Proof :}  We prove the result for $i=0$, the other case being similar. With the notation of the proof of Lemma~\ref{lem eps2-0}, one obtains
\begin{eqnarray*}
\sum_{k \in \dG_{n}}X_k (\wh{\eps}_{2k}-{\eps}_{2k})^2&=&\sum_{k \in \dG_{n}} \delta_{2k}(\wh{\bs{\theta}}^0_n-{\bs{\theta}}^0)^t
\left(\begin{array}{cc}
X_k&X_k^2\\
X_k^2&X_k^3
\end{array}\right)
(\wh{\bs{\theta}}^0_n-{\bs{\theta}}^0)\\
&=&(\bs{M}_n^0)^t(\bs{S}_{n-1}^0)^{-1}(\bs{T}_{n}^0-\bs{T}_{n-1}^0)(\bs{S}_{n-1}^0)^{-1}\bs{M}_n^0,
\end{eqnarray*}
with 
$
\bs{T}^0_{n} =\sum_{k \in \mathbb{T}_{n}}\delta_{2k}\left(
\begin{array}{cc}
X_k&X_k^2\\
X_k^2&X_k^3
\end{array}\right).
$
Proposition~\ref{prop LGN} yields the following convergence
\begin{equation*}
\lim_{n\rightarrow\infty}\ind{|\dG_n^*|>0}(\bs{\Sigma}_n^0)^{1/2}(\bs{S}_{n}^0)^{-1}(\bs{T}_{n+1}^0-\bs{T}_{n}^0)(\bs{S}_{n}^0)^{-1}(\bs{\Sigma}_n^0)^{1/2}=\bs{\Delta^0}\indnex\quad a.s.
\end{equation*}
with a new matrix $\bs{\Delta^0}$ defined by $\bs{\Delta^0}=(m-1)\bs{\Sigma}^0(\bs{S}^0)^{-2}\bs{T}^0$. Note that this new $\bs{\Delta^0}$ is again positive definite and the matrices $\bs{\Sigma}^0_n$, $\bs{S}^0_n$, $\bs{T}^0_n$, $\bs{\Sigma}^0$, $\bs{S}^0$, $\bs{T}^0$,  and $\bs{\Delta^0}$ still commute. The end of the proof is similar to that of Lemma~\ref{lem eps2-0} with the new matrix $\bs{\Delta^0}$.
\hspace{\stretch{1}}$\Box$

\begin{Lemma}\label{lem esp2-3}
Under assumptions \emph{\textbf{(H.1-5)}} and if $\kappa\geq 4$, for all $i\in\{0,1\}$,  the almost sure convergence holds
\begin{equation*}
\lim_{n\rightarrow\infty}\ind{|\dG_n^*|>0}\frac{1}{n}\sum_{k \in \dT_{n}} X_k^2(\wh{\eps}_{2k+i}-{\eps}_{2k+i})^2=q_i(2)\indnex=(m-1)tr(\bs{\Gamma}^i(\bs{S}^i)^{-2}\bs{W}^i)\indnex,
\end{equation*}
where $\bs{W}^i$ is the $2\times 2$ matrix defined by
$
\bs{W}^i=\left(
\begin{array}{cc}
\ell_i(2)&\ell_i(3)\\
\ell_i(3)&\ell_i(4)
\end{array}\right).
$
\end{Lemma}
Lemmas~\ref{lem eps2-0}, \ref{lem eps2-1}, \ref{lem esp2-2} and \ref{lem esp2-3} give the almost sure convergence of the sequence $(\bs{P}_n^{\bs{\sigma}})$.

\begin{Lemma}\label{lem lim Psigma}
Under assumptions \emph{\textbf{(H.1-5)}} and if $\kappa\geq 4$, the following convergence holds
\begin{equation*}
\lim_{n\rightarrow\infty}\ind{|\dG_n^*|>0}\frac{1}{n}\bs{P}_n^{\bs{\sigma}}=\big(q_0(0)+q_1(0), 2q_0(1),  2q_1(1), q_0(2)+q_1(2)\big)^t\indnex\quad a.s.
\end{equation*}
\end{Lemma}
It remains to give the limit of the sequence $(\bs{R}_n^{\bs{\sigma}})$.

\begin{Lemma}\label{lem lim Rsigma}
Under assumptions \emph{\textbf{(H.1-5)}} and if $\kappa\geq 8$, the following convergence holds
\begin{equation*}
\lim_{n\rightarrow\infty}\ind{|\dG_n^*|>0}\frac{1}{n}\bs{R}_n^{\bs{\sigma}}=0\qquad a.s.
\end{equation*}
\end{Lemma}

\noindent\textbf{Proof :}  It is sufficient to prove that $(\bs{R}_n^{\bs{\sigma}})$ is a martingale and that its predictable quadratic variation is almost surely $\cO(n)$. For all $k\in\dG_n$, one obtains
\begin{equation*}
\dE[{\eps}_{2k}(\wh{\eps}_{2k}-{\eps}_{2k})\ |\ \cF_n^{\cO}]=\delta_{2k}\big((a-\wh{a}_n)+(b-\wh{b}_n)X_k\big)\dE[{\eps}_{2k}\ |\ \cF_n^{\cO}]=0,
\end{equation*}
and we have the same result for the other entries of $\bs{R}_n^{\bs{\sigma}}$. Hence,  $(\bs{R}_n^{\bs{\sigma}})$ is a $(\cF_n^{\cO})$-martingale. It is also square-integrable. We are going to study $(\bs{R}_n^{\bs{\sigma}})$ component-wise.   We give the details for the last entry, the others being treated similarly. For $i\in\{0,1\}$, set
\begin{equation*}
{Q}_n^i=\sum_{\ell=1}^{n-1}(\bs{\theta}^i-\wh{\bs{\theta}}_n^i)^t\sum_{k\in\dG_{\ell}}\delta_{2k+i}\left(\begin{array}{c}X_k^2\\X_k^3\end{array}\right)\eps_{2k+i}.
\end{equation*}
The last entry of $\bs{R}_n^{\bs{\sigma}}$ can then be rewritten as ${Q}_n^0+{Q}_n^1$. The processes ${Q}_n^i$ are clearly $(\cF_n^{\cO})$-martingales with predictable quadratic variation equal to
\begin{equation*}
<{Q}^i>_n=\sum_{\ell=1}^{n-1}(\bs{M}^i_{\ell})^{t}(\bs{S}^i_{\ell-1})^{-1}\bs{\Delta}_{\ell}^i(\bs{S}^i_{\ell-1})^{-1}\bs{M}^i_{\ell},
\end{equation*}
with
$
\bs{\Delta}_n^i=\sum_{k\in\dG_n}\delta_{2k}(\sigma^2_{\veps}+2\rho_{ii}X_k+\sigma^2_{\eta}X_k^2)\left(
\begin{array}{cc}
X_k^4&X_k^5\\
X_k^5&X_k^6
\end{array}
\right).
$
Thanks to Proposition~\ref{prop LGN}, the sequence of matrices $(\bs{\Sigma}^i_n)^{1/2}(\bs{S}^i_{n-1})^{-1}\bs{\Delta}_{n}^i(\bs{S}^i_{n-1})^{-1}(\bs{\Sigma}^i_n)^{1/2}$ converges  almost surely on the non-extinction set $\overline{\cE}$ to a fixed positive definite matrix $\bs{\Delta}^i$. We now use Theorem~\ref{meta martingale} along the same lines as in the proof of Lemma~\ref{lem eps2-0} to obtain that $<{Q}^i>_n=\cO(n)$, and thus ${Q}^i_n=o(n)$. The other entries of $(\bs{R}_n^{\bs{\sigma}})$ are dealt with similarly, yielding the result.
\hspace{\stretch{1}}$\Box$\\

\noindent\textbf{Proof of  Theorem~\ref{th sigma}} It is a direct consequence of Eq.~(\ref{decomp sigma}), Proposition~\ref{prop lim SUV} and Lemmas~\ref{lem sigman}, \ref{lem lim Psigma} and \ref{lem lim Rsigma}.
\hspace{\stretch{1}}$\Box$
 
\subsection{Rate of convergence for $\bs{\wh{\rho}}_n$}
\label{section rho}
We proceed again in two steps. Recall that 
\begin{equation*}
\wh{\bs{\rho}}_n= 
\bs{V}^{-1}_{n-1}
\sum_{k \in \mathbb{T}_{n-1}}
\left( 
\wh{\eps}_{2k} \wh{\eps}_{2k+1},
2X_k \wh{\eps}_{2k}  \wh{\eps}_{2k+1},
X^2_k\wh{\eps}_{2k}  \wh{\eps}_{2k+1} 
\right)^t,
\end{equation*}
is our estimator of $\bs{\rho}=(\rho_{\veps}, \rho, \rho_{\eta})^t$, and
\begin{equation*}
{{\bs{\rho}}}_n\ =\  
\bs{V}^{-1}_{n-1}
\sum_{k \in \mathbb{T}_{n-1}}
\left( 
{\eps}_{2k} {\eps}_{2k+1},
2X_k {\eps}_{2k}  {\eps}_{2k+1},
X^2_k{\eps}_{2k}  {\eps}_{2k+1} 
\right)^t.
\end{equation*}
Our first step is to prove the convergence of ${{\bs{\rho}}}_n$ to ${{\bs{\rho}}}$. The second step is the convergence of $\wh{{\bs{\rho}}}_n-{{\bs{\rho}}}_n$ with a convergence rate.

\subsubsection{Convergence of ${{\bs{\rho}}}_n$ }
\label{section cv rho n}

The convergence of  ${{\bs{\rho}}}_n$ to ${{\bs{\rho}}}$ is again directly obtained using the standard law of large numbers for square-integrable vector-valued martingales. Note that one could also obtain a convergence rate under stronger moment assumptions using Theorem~\ref{meta martingale}.

\begin{Lemma}\label{lem rhon}
Under assumptions \emph{\textbf{(H.1-5)}} and if $\kappa\geq 8$, the following convergence holds
\begin{equation*}
\lim_{n\rightarrow\infty} \ind{|\dG_n^*|>0}{{\bs{\rho}}}_n={{\bs{\rho}}}\indnex\qquad a.s.
\end{equation*}
\end{Lemma}

\noindent\textbf{Proof :} Set 
\begin{equation*}
\bs{M}^{\bs{\rho}}_n = \bs{V}_{n-1}({{\bs{\rho}}}_n-{{\bs{\rho}}})=\sum_{\ell=1}^{n-1}\sum_{k \in \mathbb{G}_{\ell}}
\left( \begin{array}{c}
{\eps}_{2k}{\eps}_{2k+1} -\dE[{\eps}_{2k} {\eps}_{2k+1}\ |\ \cF_{\ell}^{\cO}]\\
2X_k\big({\eps}_{2k}{\eps}_{2k+1} -\dE[{\eps}_{2k} {\eps}_{2k+1}\ |\ \cF_{\ell}^{\cO}]\big)\\
X_k^2\big({\eps}_{2k}{\eps}_{2k+1} -\dE[{\eps}_{2k} {\eps}_{2k+1}\ |\ \cF_{\ell}^{\cO}]\big)
\end{array}\right).
\end{equation*}
Hence, $(\bs{M}^{\bs{\rho}}_n)$ is a square-integrable $(\cF_n^{\cO})$-martingale. One can compute all the entries of its predictable quadratic variation and prove that they all equal a constant (depending on the moments of $(\veps_2,\eta_2,\veps_3,\eta_3)$ up to order 4) multiplied by $\sum\delta_{2k}\delta_{2k+1}X_k^q$ with $q\leq 8$. Therefore, the laws of large numbers given in Proposition~\ref{prop LGN}  ensure that  $m^{-n}<\bs{M}^{\bs{\rho}}>_n$ converges almost surely to a constant matrix on the non-extinction set $\overline{\cE}$. The standard law of large numbers for square-integrable martingales then implies that  $(\bs{M}^{\bs{\rho}}_n)=o(m^n)$ a.s. 
Besides, $m^{-n}\bs{V}_n$ also converges to a fixed matrix on the non-extinction set $\overline{\cE}$ by Proposition~\ref{prop lim SUV}. Thus ${{\bs{\rho}}}_n-{{\bs{\rho}}}= \bs{V}_{n-1}^{-1}\bs{M}^{\bs{\rho}}_n$ tends to $0$ a.s. on $\overline{\cE}$ when $n$ tends to infinity.
\hspace{\stretch{1}}$\Box$

\subsubsection{Convergence of $\wh{{\bs{\rho}}}_n-{{\bs{\rho}}}_n$ }
\label{section cv rho hat n}
We now turn to the convergence of  $\wh{{\bs{\rho}}}_n-{{\bs{\rho}}}_n$. We follow the same steps as in Section~\ref{section cv sigma hat n}. One can rewrite $\bs{V}_{n-1}(\wh{{\bs{\rho}}}_n-{{\bs{\rho}}}_n)$ as
\begin{eqnarray}\label{decomp rho}
\bs{V}_{n-1}(\wh{{\bs{\rho}}}_n-{{\bs{\rho}}}_n) &= &\bs{P}_n^{\bs{\rho}}+\bs{R}_n^{\bs{\rho}},\quad \text{with}\\
\bs{P}_n^{\bs{\rho}} & = &\sum_{k \in \mathbb{T}_{n-1}}
(\wh{\eps}_{2k}-{\eps}_{2k})(\wh{\eps}_{2k+1} -{\eps}_{2k+1} )\left( 
1,
2X_k ,
X^2_k
\right)^t, \quad \text{and} \nonumber \\
\bs{R}_n^{\bs{\rho}} & = &\sum_{k \in \mathbb{T}_{n-1}}
\big(\eps_{2k+1}(\wh{\eps}_{2k}-{\eps}_{2k})+\eps_{2k}(\wh{\eps}_{2k+1} -{\eps}_{2k+1} )\big)
\left(
1,
2X_k ,
X^2_k
\right)^t. \nonumber
\end{eqnarray}
We are going to study separately the asymptotic properties of $\bs{P}_n^{\bs{\rho}}$ and $\bs{R}_n^{\bs{\rho}}$. The limit of $\bs{R}_n^{\bs{\rho}}$ is obtained as in Lemma~\ref{lem lim Rsigma}.

\begin{Lemma}\label{lem lim Rrho}
Under assumptions \emph{\textbf{(H.1-5)}} and if $\kappa\geq 8$, one obtains
\begin{equation*}
\lim_{n\rightarrow\infty}\ind{|\dG_n^*|>0}\frac{1}{n}\bs{R}_n^{\bs{\rho}}=0\qquad a.s.
\end{equation*}
\end{Lemma}

\begin{Lemma}\label{lem P0}
Under assumptions \emph{\textbf{(H.1-5)}} and if $\kappa\geq 4$,  the following almost sure convergence holds
\begin{equation*}
\lim_{n\rightarrow\infty}\frac{\ind{|\dG_n^*|>0}}{n}\sum_{k \in \dT_{n}} (\wh{\eps}_{2k}-{\eps}_{2k})(\wh{\eps}_{2k+1}-{\eps}_{2k+1})=q_{01}(0)\indnex=\frac{m-1}{2}tr(\bs{\Gamma}\bs{S}^{-2}\bs{J}^{01})\indnex
\end{equation*}
where 
$
\bs{J}^{01}=\left(\begin{array}{cc}
0&\bs{S}^{01}\\
\bs{S}^{01}&0
\end{array}\right)$, and
$\bs{S}^{01}=\left(
\begin{array}{cc}
\ell_{01}(0)&\ell_{01}(1)\\
\ell_{01}(1)&\ell_{01}(2)
\end{array}\right)$.
\end{Lemma}

\noindent\textbf{Proof :} First, notice that for all $k\in\dG_n$, the following decomposition holds
$$2 (\wh{\eps}_{2k}-{\eps}_{2k})(\wh{\eps}_{2k+1}-{\eps}_{2k+1})
 = \delta_{2k}\delta_{2k+1}(\wh{\bs{\theta}}_n-\bs{\theta})^t
\left(\begin{array}{cccc}
0&0&1&X_k\\
0&0&X_k&X_k^2\\
1&X_k&0&0\\
X_k&X_k^2&0&0
\end{array}\right)(\wh{\bs{\theta}}_n-\bs{\theta}).$$
Hence, one obtains
\begin{equation*}
2\sum_{k \in \dT_{n}} (\wh{\eps}_{2k}-{\eps}_{2k})(\wh{\eps}_{2k+1}-{\eps}_{2k+1})=\sum_{\ell=1}^n\bs{M}_{\ell}^t\bs{S}_{\ell-1}^{-1}(\bs{J}^{01}_{\ell}-\bs{J}^{01}_{\ell-1})\bs{S}_{\ell-1}^{-1}\bs{M}_{\ell},
\end{equation*}
with 
$
\bs{J}^{01}_{n}=\left(\begin{array}{cc}
0&\bs{S}^{01}_n\\
\bs{S}^{01}_n&0
\end{array}\right)$
 and $\bs{S}^{01}_n=\sum_{k\in\dT_n}2\delta_{2k}\delta_{2k+1}
\left(\begin{array}{cc}
1&X_k\\
X_k&X_k^2
\end{array}\right)
$.
Set $\bs{\Delta}_n=\bs{\Sigma}_{n}^{1/2}\bs{S}_{n}^{-1}(\bs{J}^{01}_{n+1}-\bs{J}^{01}_{n})\bs{S}_{n}^{-1}\bs{\Sigma}_{n}^{1/2}$. Proposition~\ref{prop LGN} yields
\begin{equation*}
\lim_{n\rightarrow\infty}\ind{|\dG_n^*|>0}\bs{\Delta}_n=\bs{\Delta}\indnex=(m-1)\bs{\Sigma}^{1/2}\bs{S}^{-1}\bs{J}^{01}\bs{S}^{-1}\bs{\Sigma}^{1/2}\indnex\quad a.s.
\end{equation*} 
Hence, as in the proof of Lemma~\ref{lem eps2-0}, it is sufficient to study the convergence of
$
\sum_{\ell=1}^n\bs{M}_{\ell}^t\bs{\Sigma}_{\ell-1}^{-1/2}\bs{\Delta}\bs{\Sigma}_{\ell-1}^{-1/2}\bs{M}_{\ell}.
$
To this end, we apply Theorem~\ref{meta martingale} to the martingale $(\bs{M}_n)$ and with the sequence of matrices $\bs{\Xi}_n=\bs{\Sigma}_{n}^{1/2}\bs{\Delta}^{-1}\bs{\Sigma}_{n}^{1/2}$. As we have seen before, assumptions \textbf{(A.1-3)} and \textbf{(A.6)} hold. Assumption~\textbf{(A.4)} is a direct consequence of Lemma~\ref{lem lim Sigma}. We will not investigate Assumption~\textbf{(A.5)} but directly prove Lemma \ref{lemlimB}, a part of Theorem~\ref{meta martingale},  in our specific context that is
$$ \cB'_{n+1}=2\sum_{k=1}^n \bs{M}_k^t\bs{\Xi}_k^{-1}\Delta \bs{M}_{k+1} = o(n).$$
Note that the matrix $\bs{J}^{01}$ has a special property, namely
\begin{equation*}
(\bs{J}^{01})^2=\left(\begin{array}{cc}
(\bs{S}^{01})^2&0\\
0&(\bs{S}^{01})^2
\end{array}\right)=
\left(\begin{array}{cc}
\bs{S}^{01}&0\\
0&\bs{S}^{01}
\end{array}\right)^2=(\bs{I}^{01})^2,
\end{equation*}
so that although $\bs{J}^{01}$ is not positive definite, $\bs{I}^{01}$ is. As a result, as the matrices $\bs{J}^{01}_n$, $\bs{\Sigma}_n$ and $\bs{S}_n$ commute, one has $\bs{\Delta}^2=(\bs{\Delta}')^2$ with the positive definite matrix
$
\bs{\Delta}'=(m-1)\bs{\Sigma}^{1/2}\bs{S}^{-1}\bs{I}^{01}\bs{S}^{-1}\bs{\Sigma}^{1/2}
$.
We have seen that the coefficient $\alpha$ given by Lemma~\ref{lem A5} satisfies
$
\bs{\Sigma}_{n}^{-1}\bs{\Gamma}_n\bs{\Sigma}_{n}^{-1}\leq\alpha\big(\bs{\Sigma}_{n-1}^{-1}-\bs{\Sigma}_{n}^{-1}\big).
$
In view of the property of $\bs{J}^{01}$ and the fact that the matrices $\bs{\Delta}$, $\bs{\Delta}'$, $\bs{\Sigma}_n$ and $\bs{\Gamma}_n$ commute, we obtain
\begin{eqnarray*}
\big(\bs{\Sigma}_{n}^{1/2}\bs{\Delta}^{-1}\bs{\Sigma}_{n}^{1/2}\big)^{-1}\bs{\Gamma}_n\big(\bs{\Sigma}_{n}^{1/2}\bs{\Delta}^{-1}\bs{\Sigma}_{n}^{1/2}\big)^{-1}
&=&\bs{\Delta}'\bs{\Sigma}_{n}^{-1}\bs{\Gamma}_n\bs{\Sigma}_{n}^{-1}\bs{\Delta}'\\
&\leq&\alpha\bs{\Delta}'\big(\bs{\Sigma}_{n-1}^{-1}-\bs{\Sigma}_{n}^{-1}\big)\bs{\Delta}'.
\end{eqnarray*}
Thanks to Lemma \ref{lem lim Sigma}, we have the convergence
$$\lim_{n \rightarrow \infty} \frac{\lambda_{max}(\bs{\Sigma}_{n-1}^{-1}-\bs{\Sigma}_{n}^{-1}\big)}{\lambda_{min}(\bs{\Sigma}_{n-1}^{-1}-\bs{\Sigma}_{n}^{-1})} =  \frac{\lambda_{max}(\bs{\Sigma}^{-1})}{\lambda_{min}(\bs{\Sigma}^{-1})} >0.$$So 
there exists $n_0>0$ and $\beta>0$ such that 
$$\beta \lambda_{min}(\bs{\Sigma}_{n-1}^{-1}-\bs{\Sigma}_{n}^{-1}) \geq \lambda_{max}(\bs{\Sigma}_{n-1}^{-1}-\bs{\Sigma}_{n}^{-1}) \lambda_{max}((\bs{\Delta}')^2) \text{ for } n \geq n_0,$$ which implies
$\bs{\Delta}'\big(\bs{\Sigma}_{n-1}^{-1}-\bs{\Sigma}_{n}^{-1}\big)\bs{\Delta}' \leq \beta \big(\bs{\Sigma}_{n-1}^{-1}-\bs{\Sigma}_{n}^{-1}\big) $, and 
\begin{eqnarray*}
\dE[\Delta \cB_{n+1}^{'2}|\cF_n^{\cO}]=4\bs{M}_n^t\bs{\Xi}_n^{-1}\bs{\Gamma}_n\bs{\Xi}_n^{-1}\bs{M}_n \leq \alpha \beta \bs{M}_n^t \big(\bs{\Sigma}_{n-1}^{-1}-\bs{\Sigma}_{n}^{-1}\big) \bs{M}_n \text{ for } n \geq n_0.
\end{eqnarray*}$ \bs{M}_n^t \big(\bs{\Sigma}_{n-1}^{-1}-\bs{\Sigma}_{n}^{-1}\big) \bs{M}_n$ is the increment of $\cA_n$ defined in Eq \ref{maindecomart} when $\Xi_n = \Sigma_n$. By the arguments used in the proof of Proposition \ref{prop vitesse M}, 
$\sum_{\ell = 1}^n\bs{M}_\ell^t \big(\bs{\Sigma}_{\ell-1}^{-1}-\bs{\Sigma}_{\ell}^{-1}\big) \bs{M}_\ell = \cO(n)$ which implies $ \cB'_{n+1}=o(n)$. We can then apply Theorem~\ref{meta martingale} which immediately yields the result.
\hspace{\stretch{1}}$\Box$\\

\noindent Similarly, one obtains the two following convergences.

\begin{Lemma}\label{lem P1}
Under assumptions \emph{\textbf{(H.1-5)}} and if $\kappa\geq 4$,  the following  almost sure convergence holds
\begin{equation*}
\lim_{n\rightarrow\infty}\frac{\ind{|\dG_n^*|>0}}{n}\sum_{k \in \dT_{n}} X_k(\wh{\eps}_{2k}-{\eps}_{2k})(\wh{\eps}_{2k+1}-{\eps}_{2k+1})=q_{01}(1)\indnex=\frac{m-1}{2}tr(\bs{\Gamma}\bs{S}^{-2}\bs{K}^{01})\indnex
\end{equation*}
where 
$
\bs{K}^{01}=\left(\begin{array}{cc}
0&\bs{T}^{01}\\
\bs{T}^{01}&0
\end{array}\right)$ and
$\bs{T}^{01}=\left(
\begin{array}{cc}
\ell_{01}(1)&\ell_{01}(2)\\
\ell_{01}(2)&\ell_{01}(3)
\end{array}\right).
$
\end{Lemma}

\begin{Lemma}\label{lem P2}
Under assumptions \emph{\textbf{(H.1-5)}} and if $\kappa\geq 4$,  the following  almost sure convergence holds
\begin{equation*}
\lim_{n\rightarrow\infty}\frac{\ind{|\dG_n^*|>0}1}{n}\sum_{k \in \dT_{n}} X_k^2(\wh{\eps}_{2k}-{\eps}_{2k})(\wh{\eps}_{2k+1}-{\eps}_{2k+1})=q_{01}(2)\indnex=\frac{m-1}{2}tr(\bs{\Gamma}\bs{S}^{-2}\bs{L}^{01})\indnex
\end{equation*}
where 
$
\bs{L}^{01}=\left(\begin{array}{cc}
0&\bs{W}^{01}\\
\bs{W}^{01}&0
\end{array}\right)$, and
$\bs{W}^{01}=\left(
\begin{array}{cc}
\ell_{01}(2)&\ell_{01}(3)\\
\ell_{01}(3)&\ell_{01}(4)
\end{array}\right).
$
\end{Lemma}
Thus, one obtains the following limit.

\begin{Lemma}\label{lem lim Prho}
Under assumptions \emph{\textbf{(H.1-5)}} and if $\kappa\geq 4$,  the following  almost sure convergence holds
\begin{equation*}
\lim_{n\rightarrow\infty}\ind{|\dG_n^*|>0}\frac{1}{n}\bs{P}_n^{\bs{\sigma}}=\big(q_{01}(0),  2q_{01}(1), q_{01}(2)\big)^t\indnex\qquad a.s.
\end{equation*}
\end{Lemma}

\noindent\textbf{Proof of  Theorem~\ref{th rho}} It is a direct consequence of Eq.~(\ref{decomp rho}), Proposition~\ref{prop lim SUV} and Lemmas~\ref{lem rhon}, \ref{lem lim Rrho} and \ref{lem lim Prho}.
\hspace{\stretch{1}}$\Box$

\section{Asymptotic normality}
\label{section TCL}
To derive the central limit theorems (CLT), we use a CLT for martingales given in \cite[Theorem~3.II.10]{Duflo97}. To this aim, we use a new filtration. Namely, instead of using the observed generation-wise filtration, we will use the sister pair-wise one. Let
\begin{equation*}
\mathcal{G}^{\cO}_p=\cO \vee \sigma\{\delta_1X_1,\ (\delta_{2k}X_{2k}, \delta_{2k+1}X_{2k+1}),\ 1\leq k\leq p\}
\end{equation*}
be the $\sigma$-algebra generated by the whole history $\cO$ of the Galton-Watson process and all observed individuals up to the offspring of individual $p$. Hence $(\delta_{2k}\eps_{2k}, \delta_{2k+1}\eps_{2k+1})$ is $\mathcal{G}^{\cO}_k$-measurable.
In all the sequel, we will work on the non-extinction probability space $(\overline{\cE},\mathbb{P}_{\overline{\cE}})$ and we denote by $\mathbb{E}_{\overline{\cE}}$ the corresponding expectation.\\

\noindent\textbf{Proof of Theorem~\ref{thmCLT}, first step: } For a fixed integer $n \geq 1$, let us define
 the $\mathcal{G}^{\cO}_p$-martingale $(\bs{M}^{(n)}_p)_{\{p \geq 1\}} $ by 
\begin{equation*}
\bs{M}^{(n)}_p = \frac{1}{{|\dT^*_n|}^{1/2}} \sum_{k=1}^{p}\bs{D}_k
\qquad \text{with} \qquad 
\bs{D}_k=\left( 
\eps_{2k},
X_k\eps_{2k},
\eps_{2k+1},
X_k\eps_{2k+1}
\right)^t.
\end{equation*}
Indeed, under {(\textbf{H.1-5})}, $D_k$ is clearly a $\mathcal{G}^{\cO}_k$-martingale difference sequence. 
Set $\nu_n=|\dT_n|=2^{n+1}-1$. Therefore the following equality holds
\begin{equation*}
\bs{M}^{(n)}_{\nu_n}=  \frac{1}{{|\dT^*_n|}^{1/2}} \sum_{k=1}^{|\dT_n|}\bs{D}_k= \frac{1}{{|\dT^*_n|}^{1/2}} \bs{M}_n.
\end{equation*}
We want to apply Theorem~3.II.10 of \cite{Duflo97} to the process $(\bs{M}^{(n)}_{\nu_n})$.
As the non-extinction set $\overline{\cE}$ is in $\mathcal{G}^{\cO}_{k}$ for every $k \geq 1$, it is easy to prove that 
\begin{equation*}
\mathbb{E}_{\overline{\cE}}[\bs{D}_k\bs{D}_k^t|\mathcal{G}^{\cO}_{k-1}]=\mathbb{E}[\bs{D}_k\bs{D}_k^t|\mathcal{G}^{\cO}_{k-1}] =\bs{\gamma}_k
\otimes\left(
\begin{array}{cc}
1&X_k\\
X_k&X_k^2
\end{array}
\right),
\end{equation*}where $\gamma_k$ is defined in Eq. (\ref{gamk}). 
Lemma \ref{lem lim Gamma} gives  the $\mathbb{P}_{\overline{\cE}}$ almost sure limit of the following process
\begin{equation*}
\frac{1}{|\dT^*_n|} \sum_{k \in \dT^*_n}\mathbb{E}_{\overline{\cE}}[\bs{D}_k\bs{D}^t_k|\mathcal{G}^{\cO}_{k-1}]    \xrightarrow[n\rightarrow\infty]{} \bs{\Gamma} \hspace{1 cm} \text{a.s.}
\end{equation*}
Therefore, the first assumption of Theorem~3.II.10 of \cite{Duflo97} holds under
$\mathbb{P}_{\overline{\cE}}$.  Thanks to $\gamma \geq \kappa \geq 4$, we can easily  prove that for some $r>2$, one obtains
\begin{equation*}
\sup_{k\geq 0}
\dE[\|\bs{D}_k\|^r|\mathcal{G}^{\cO}_{k-1}]<\infty
\hspace{1cm}\text{a.s.}
\end{equation*} 
which in turn implies the Lindeberg condition. We can now conclude that under
$\mathbb{P}_{\overline{\cE}}$ the following convergence holds
\begin{equation*}
\frac{1}{|\dT^*_{n-1}|^{1/2}}\sum_{k \in \mathbb{T}_{n-1}}\bs{D}_k=\frac{1}{|\dT^*_{n-1}|^{1/2}}\bs{M}_n
\liml
\cN(0,\bs{\Gamma}).
\end{equation*}
Finally, result (\ref{CLTtheta}) follows from
Eq.~(\ref{diff}) and Proposition~\ref{prop lim SUV} together with Slutsky's Lemma.
\hspace{\stretch{1}}$\Box$\ \\

\noindent\textbf{Proof of Theorem~\ref{thmCLT}, second step:} We apply Theorem~3.II.10 of \cite{Duflo97} again
to the sequences  $(\bs{M}^{(\bs{\sigma},n)}_p)_{\{p \geq 1\}} $ of $\mathcal{G}^{\cO}_p$-martingales defined by 
\begin{equation*}
|\dT^*_{n-1}|^{1/2}\bs{M}^{(\bs{\sigma}, n)}_{p} 
=  \sum_{k =1}^p  \bs{D}_k^{\bs{\sigma}}\\
 =   \sum_{k =1}^p
\left( \begin{array}{c}
{\eps}^2_{2k} +{\eps}^2_{2k+1} -\dE[{\eps}^2_{2k} +{\eps}^2_{2k+1}\ |\ \cF_{r_k}^{\cO}]\\
2X_k({\eps}^2_{2k} -\dE[{\eps}^2_{2k}\ |\ \cF_{r_k}^{\cO}])\\
2X_k({\eps}^2_{2k+1}-\dE[{\eps}^2_{2k+1}\ |\ \cF_{r_k}^{\cO}])\\
X^2_k( {\eps}^2_{2k} + {\eps}^2_{2k+1}-\dE[{\eps}^2_{2k} +{\eps}^2_{2k+1}\ |\ \cF_{r_k}^{\cO}] )
\end{array}\right),
\end{equation*}
where $r_k$ denotes the generation of $k$. Set $\nu_n=|\dT_{n-1}|=2^{n}-1$. One obtains
$
{|\dT^*_{n-1}|}^{1/2} \bs{M}^{(\bs{\sigma}, n)}_{\nu_n} 
 = \bs{U}_{n-1}({{\bs{\sigma}}}_n-{{\bs{\sigma}}})
$.
We have to study the limit $\bs{\Gamma}^{\sigma}$ of the process 
$\frac{1}{|\dT^*_{n-1}|} \sum_{k =1}^{\dT_{n-1}}  \mathbb{E}_{\overline{\cE}}[\bs{D}_k^{\bs{\sigma}}(\bs{D}_k^{\bs{\sigma}})^t\ |\ \mathcal{G}^{\cO}_{k-1}].$ In order to compute the conditional expectation, let us denote, for $k\geq 1$
\begin{eqnarray*}
A_i(k)
& = &  \delta_{2k+i} \Big( \sum_{r=0}^4 C_4^l \theta((1-i)(4-r),(1-i)r,i(4-r),  ir)X_k^r  \\
&  &    - \big(\sigma^4_\veps + 4 \rho_{ii}\sigma^2_\veps X_k+(4 \rho_{ii}^2+2 \sigma^2_\veps \sigma^2_\eta)X_k^2 +4 \rho_{ii}\sigma^2_\eta X_k^3+ \sigma^4_\eta X_k^4\big)\Big), \\
A_{01}(k) & =  &   \delta_{2k}\delta_{2k+1} \Big( \sum_{r=0}^2 \sum_{s=0}^2 C_2^r C_2^{s} \theta(2-r,r,2-s,s)X_k^{r+s}\\
&  & - \big( \sigma^4_\veps + 2 \sigma^2_\veps (\rho_{00}+ \rho_{11})X_k  +(2 \sigma^2_\veps \sigma^2_\eta +4 \rho_{00}\rho_{11})X_k^2+2 \sigma^2_\eta (\rho_{00}+ \rho_{11})X^3_k +  \sigma^4_\eta X_k^4 \big)\Big ),
\end{eqnarray*}
and $B_i(k)=A_i(k)+A_{01}(k)$. Using these notations,  we obtain
\begin{eqnarray*}
\lefteqn{\mathbb{E}_{\overline{\cE}}[\bs{D}_k^{\bs{\sigma}}(\bs{D}_k^{\bs{\sigma}})^t\ |\ \mathcal{G}^{\cO}_{k-1}]}\\
&=& \left( \begin{array}{cccc}
(B_0 + B_1)(k)  & 2X_kB_0(k)  & 2X_kB_1(k)   & X_k^2 (B_0+ B_1)(k)\\
2X_kB_0(k)  & 4 X_k^2A_0(k)  & 4 X_k^2A_{01}(k) &  2X^3_kB_0(k) \\
2X_kB_1(k) &  4 X_k^2A_{01}(k) & 4 X_k^2A_1(k)   & 2X^3_kB_1(k) \\
X_k^2 (B_0+ B_1)(k) & 2X^3_kB_0(k)  & 2X^3_kB_1(k) & X_k^4(B_0 + B_1)(k) 
\end{array}\right)
\end{eqnarray*}
The almost sure limit of the above quantity is given by Proposition \ref{prop LGN}. Indeed, the following convergences hold 
\begin{equation*}
\lim_{n\rightarrow\infty}\frac{1}{|\dT^*_{n-1}|} \sum_{k =1}^{\dT_{n-1}} A_i(k)X^q_k = A_i^q \quad \text{and}  \quad \lim_{n\rightarrow\infty}\frac{1}{|\dT^*_{n-1}|} \sum_{k =1}^{\dT_{n-1}} A_{01}(k)X^q_k = A_{01}^q \quad \text{a.s.}
\end{equation*}
with
\begin{eqnarray*}
A_i^q
& = &  \sum_{r=0}^4 C_4^l \theta((1-i)(4-r),(1-i)r,i(4-r),  ir)\ell_i(r+q)  \\
&  &    - \big(\sigma^4_\veps\ell_i(q) + 4 \rho_{ii}\sigma^2_\veps \ell_i(1+q)+(4 \rho_{ii}^2+2 \sigma^2_\veps \sigma^2_\eta)\ell_i(2+q) +4 \rho_{ii}\sigma^2_\eta \ell_i(3+q)+ \sigma^4_\eta \ell_i(4+q)4\big), \\
A_{01}^q & =  &  \sum_{r=0}^2 \sum_{s=0}^2 C_2^r C_2^{s} \theta(2-r,r,2-s,s)\ell_{01}(r+s+q)\\
&  & - \big( \sigma^4_\veps\ell_{01}(q) + 2 \sigma^2_\veps (\rho_{00}+ \rho_{11})\ell_{01}(1+q)  +(2 \sigma^2_\veps \sigma^2_\eta +4 \rho_{00}\rho_{11})\ell_{01}(2+q)\\
&& +2 \sigma^2_\eta (\rho_{00}+ \rho_{11})\ell_{01}(3+q) +  \sigma^4_\eta \ell_{01}(4+q) \big).
\end{eqnarray*}
We also set $B_i^q=A_i^q+A_{01}^q$. With these notations, we are able to explicit the limit matrix $\bs{\Gamma}^{\sigma}$ of the process 
$\frac{1}{|\dT^*_{n-1}|} \sum_{k =1}^{\dT_{n-1}}  \mathbb{E}_{\overline{\cE}}[\bs{D}_k^{\bs{\sigma}}(\bs{D}_k^{\bs{\sigma}})^t\ |\ \mathcal{G}^{\cO}_{k-1}]$
\begin{equation}
\label{Gammasigma}
\bs{\Gamma}^{\sigma} = \left( 
\begin{array}{cccc}
B^0_0 + B^0_1 & 2B^1_0       & 2B^1_1       & B^2_0+ B^2_1\\
2B^1_0               & 4 A^2_0      & 4 A^2_{01} &  2B^3_0 \\
2B^1_1              &  4 A^2_{01} & 4 A_1^2     & 2B^3_1\\
B^2_0+ B^2_1 & 2B^3_0        & 2B^3_1       & B^4_0 + B^4_1 
\end{array}
\right)
\end{equation}The first assumption of Theorem~3.II.10 of \cite{Duflo97} holds under
$\mathbb{P}_{\overline{\cE}}$ and we easily prove the second one to conclude 
that under $\mathbb{P}_{\overline{\cE}}$ 
\begin{equation*}
\bs{M}^{\bs{\sigma},n}_{\nu_n}={{|\dT^*_{n-1}|^{-1/2}}}\sum_{k \in \mathbb{T}_{n-1}}\bs{D}_k^{\bs{\sigma}}
={|\dT^*_{n-1}|^{-1/2}} \bs{U}_{n-1}({{\bs{\sigma}}}_n-{{\bs{\sigma}}})
\liml
\cN(0,\bs{\Gamma}^\sigma).
\end{equation*}
We conclude using Proposition  \ref{prop lim SUV} and Theorem \ref{th sigma}.
\hspace{\stretch{1}}$\Box$\ \\

\noindent\textbf{Proof of Theorem~\ref{thmCLT}, third step:} We apply again Theorem~3.II.10 of \cite{Duflo97} with
to the sequence of $\mathcal{G}^{\cO}_p$-martingales $(\bs{M}^{(\bs{\rho} ,n)}_p)_{\{p \geq 1\}} $ defined by 
\begin{equation*}
|\dT^*_{n-1}|^{1/2}\bs{M}^{(\bs{\rho}, n)}_{p} = \sum_{k =1}^p  \bs{D}_k^{\bs{\rho}}
=\sum_{k = 1}^p
\left( \begin{array}{c}
{\eps}_{2k}{\eps}_{2k+1} -\dE[{\eps}_{2k} {\eps}_{2k+1}\ |\ \cF_{r_k}^{\cO}]\\
2X_k\big({\eps}_{2k}{\eps}_{2k+1} -\dE[{\eps}_{2k} {\eps}_{2k+1}\ |\ \cF_{r_k}^{\cO}]\big)\\
X_k^2\big({\eps}_{2k}{\eps}_{2k+1} -\dE[{\eps}_{2k} {\eps}_{2k+1}\ |\ \cF_{r_k}^{\cO}]\big)
\end{array}\right).
\end{equation*}
Set $\nu_n=|\dT_{n-1}|=2^{n}-1$. Thus one can rewrite 
$
{|\dT^*_{n-1}|^{1/2}} \bs{M}^{(\rho, n)}_{\nu_n} 
 = \bs{V}_{n-1}({{\bs{\rho}}}_n-{{\bs{\rho}}})
$.
Let us denote
\begin{eqnarray*}
C(k)
& = &    \delta_{2k}\delta_{2k+1} \Big( (\theta(2,0,2,0)-\rho_{\veps}^2) + 2\big(\theta(2,0,1,1)+ \theta(1,1,2,0)- 2 \rho \rho_\veps\big)X_k\\
 &  &   +  \big(\theta(0,2,2,0) + \theta(2,0,0,2) + 4 \theta(1,1,1,1) - 4 \rho^2-2 \rho_\veps\ \rho_\eta\big)X_k^2\\
&  &  +2  \big(\theta(0,2,1,1) + \theta(1,1,0,2)- 2 \rho \rho_\eta\big)X^3_k + \big(\theta(2,0,2,0)-\rho_\eta^2 \big)X_k^4\Big).
\end{eqnarray*}
We are now able to write
\begin{equation*}
\mathbb{E}_{\overline{\cE}}[\bs{D}_k^{\bs{\rho}}(\bs{D}_k^{\bs{\rho}})^t\ |\ \mathcal{G}^{\cO}_{k-1}] = C(k)
\left( 
\begin{array}{ccc}
1  & 2X_k  &  X_k^2 \\
2X_k & 4 X_k^2  &  2X^3_k \\
X_k^2  & 2X^3_k & X_k^4
\end{array}
\right)
\end{equation*}For the determination of  the limit $\bs{\Gamma}^{\bs{\rho}}$ of 
${|\dT^*_{n-1}|^{-1}} \sum_{k =1}^{\dT_{n-1}}  \mathbb{E}_{\overline{\cE}}[\bs{D}_k^{\bs{\rho}}(\bs{D}_k^{\bs{\rho}})^t\ |\ \mathcal{G}^{\cO}_{k-1}]$,
 let us remark, using Proposition \ref{prop LGN} that 
\begin{equation*}
\lim_{n\rightarrow\infty}\frac{1}{|\dT^*_{n-1}|} \sum_{k =1}^{\dT_{n-1}} C(k)X^q_k = C^q  \quad \text{a.s.}
\end{equation*}with
\begin{eqnarray*}
C^q
& = & (\theta(2,0,2,0)-\rho_{\veps}^2)\ell_{01}(q) + 2\big(\theta(2,0,1,1)+ \theta(1,1,2,0)- 2 \rho \rho_\veps\big)\ell_{01}(1+q)\\
 &  &   +  \big(\theta(0,2,2,0) + \theta(2,0,0,2) + 4 \theta(1,1,1,1) - 4 \rho^2-2 \rho_\veps\ \rho_\eta\big)\ell_{01}(2+q)\\
&  &  +2  \big(\theta(0,2,1,1) + \theta(1,1,0,2)- 2 \rho \rho_\eta\big)\ell_{01}(3+q) + \big(\theta(2,0,2,0)-\rho_\eta^2 \big)\ell_{01}(4+q).
\end{eqnarray*}
The matrix $\bs{\Gamma}^{\bs{\rho}} $ is thus given by
\begin{equation}
\label{Gammarho}
\bs{\Gamma}^{\rho} = \left( 
\begin{array}{ccc}
C^0 & 2C^1  &  C^2 \\
2C^1  & 4 C^2  &  2C^3 \\
C^2  & 2C^3& C^4
\end{array}
\right)\end{equation}The first assumption of Theorem~3.II.10 of \cite{Duflo97} holds under
$\mathbb{P}_{\overline{\cE}}$ and we easily prove the second one to conclude 
that under $\mathbb{P}_{\overline{\cE}}$ 
\begin{equation*}
\bs{M}^{\bs{\rho},n}_{\nu_n}=|\dT^*_{n-1}|^{-1/2}\sum_{k \in \mathbb{T}_{n-1}}\bs{D}_k^{\bs{\rho}}=
{|\dT^*_{n-1}|^{-1/2}}\bs{V}_{n-1}({{\bs{\rho}}}_n-{{\bs{\rho}}})
\liml
\cN(0,\bs{\Gamma}^{\bs{\rho}}).
\end{equation*}
We conclude using Proposition  \ref{prop lim SUV} and Theorem \ref{th rho}.
\hspace{\stretch{1}}$\Box$

\section{Application to real data}
\label{section data}
We have applied our estimation procedure to the Escherichia coli data of  \cite{SMPT05} {(these data are available on request {from} the corresponding author of  \cite{SMPT05})}. E. coli  is a
rod-shaped bacterium that reproduces by dividing in the middle. Each
cell has thus a new end (or \emph{pole}), and an older one. The cell
that inherits the old pole of its mother is called the old pole cell,
the cell that inherits the new pole of its mother is called the new
pole cell. Therefore, each cell has a \emph{type}: old pole (even) or
new pole (odd), inducing asymmetry in the cell division. Stewart et al. \cite{SMPT05} filmed colonies of dividing cells, determining the
complete lineages and the growth rate of each cell. 
Several attempts have already been made to fit BAR processes to these data, see \cite{GBPSDT05, Guy07, DM08,SGM11}, but only with fixed coefficients models. In particular, \cite{GBPSDT05} suggests that such models cannot explain all the randomness of the data.

We have run our estimators on the data set
\texttt{penna-2002-10-04-4} from the experiments of \cite{SMPT05}.
It is the largest data set of the experiment. It contains 663 cells up
to generation 9 (note that there would be 1023 cells in a full tree up
to generation 9). {For each of the 663 observed cells, the measure of the growth rate is available. For each cell, its length was recorded from birth to division and the corresponding growth curve was fitted by an exponential function $t\mapsto\exp(\lambda t)$  where $\lambda$ is called the growth rate of the cell. Growth rates go from the minimum value $0.009$ to the maximal $0.067$. The $0.01$-quantile equals $0.024$ and the $0.99$-quantile equals $0.049$. Mean and median equal $0.037$ (std: $0.004$).}
\begin{table}[htdp]
\centering
\begin{tabular}{|c|c|c|c|}
\hline
 $a$ &$b$ &$c$ &$d$ \\
\hline
 $0.0363$ &$0.0266$& $0.0306$&$0.1706$ \\
 $[0.0275, 0.0450]$& $[-0.2094, 0.2627]$&$[0.0216, 0.0396]$ & $[-0.0709, 0.4120]$\\
\hline
\end{tabular}
\caption{Estimation of $\bs{\theta}$ on the data set \emph{penna-2002-10-04-4} }
\label{tab:1}
\end{table}
{Note that even though the number of observed generations $n=9$ is low, the rate of convergence of  our estimators is $\vert \dT^*_n \vert ^{-1/2}$ which is of order $\pi^{-n/2}$. Here for $n=9$, $\vert \dT^*_9 \vert =663$.}
Table \ref{tab:1} gives the estimation $\wh{\bs{\theta}}_9$ of ${\bs{\theta}}$
with the 95\% Confidence Interval (C.I.) of each coefficient. Note that our estimator $\wh{\bs{\theta}}_n$ of ${\bs{\theta}}$ is exactly the same as in \cite{SGM11}, and of course we obtain the same point estimation. The confidence intervals are wider, as the variance is different.
More precisely, the variance is given by the CLT for $\bs{\theta}$ in Eq.~(\ref{CLTtheta}). We have approximated it by $\bs{S}^{-1}_8\bs{\Gamma}_8 \bs{S}^{-1}_8$
thanks to the convergences given in Proposition~\ref{prop lim SUV} and Lemma~\ref{lem lim Gamma}.
\begin{table}[htdp]
\centering
\begin{tabular}{|c|c|}
\hline
$\sigma^2_{\veps}$  &$\sigma^2_{\eta}$\\
\hline
$0.0004$&$0.2431$\\
$[-0,0002, 0.0010]$&$[-0.0750, 0.5613]$\\
\hline
\end{tabular}
\caption{Estimation of noise variances on the data set \emph{penna-2002-10-04-4} }\label{tab:2}
\end{table}
Table \ref{tab:2} gives the estimation of the variance coefficients $\sigma_{\veps}^2$ and $\sigma_{\eta}^2$ of ${\bs{\sigma}}$ (other covariance coefficients of $\wh{\bs{\sigma}}_9$ and $\wh{\bs{\rho}}_9$ can be computed but are less easy to interpret). The variance of these parameters is again given by the central limit Theorem~\ref{thmCLT}. To obtain confidence intervals, one needs an estimation of the joint moments of $(\veps_2,\eta_2,\veps_3,\eta_3)$ up to the order $4$. Such estimators can be easily derived following the same ideas as in Section~\ref{section estimators}. {A Wald's test for the positivity of $\sigma_{\veps}^2$ (resp. $\sigma_{\eta}^2$) can be derived from Theorem~\ref{thmCLT}. It rejects the null hypothesis $H_0 : \sigma_{\veps}^2=0$ (resp. $H_0 : \sigma_{\eta}^2=0$) with p-value  $p=0.0799$ (resp. $p=0.0671$).} We are not far from supporting the validity of the random coefficients model.

\end{document}